\documentclass[a4paper]{amsart}

\usepackage{a4}
\usepackage[utf8]{inputenc}
\usepackage[english]{babel}

\usepackage{amsmath,amssymb,latexsym}
\usepackage{hyperref,todonotes}

\newcommand{\SP}{\vspace{0.2cm}\par}

\newtheorem{theorem}{Theorem}[section]
\newtheorem{lemma}[theorem]{Lemma}
\newtheorem{definition}[theorem]{Definition}

\newtheorem{corollary}[theorem]{Corollary}

\newtheorem{proposition}[theorem]{Proposition}

\allowdisplaybreaks 
\usepackage[numbers]{natbib} 

\newcommand\Tstrut{\rule{0pt}{2.7ex}}         
\newcommand\Bstrut{\rule[-1.3ex]{0pt}{0pt}}   

\makeatletter
\let\save@mathaccent\mathaccent
\newcommand*\if@single[3]{
	\setbox0\hbox{${\mathaccent"0362{#1}}^H$}%
	\setbox2\hbox{${\mathaccent"0362{\kern0pt#1}}^H$}%
	\ifdim\ht0=\ht2 #3\else #2\fi }
\newcommand*\rel@kern[1]{\kern#1\dimexpr\macc@kerna}
\newcommand*\widebar[1]{\@ifnextchar^{{\wide@bar{#1}{0}}}{\wide@bar{#1}{1}}}
\newcommand*\wide@bar[2]{\if@single{#1}{\wide@bar@{#1}{#2}{1}}{\wide@bar@{#1}{#2}{2}}}
\newcommand*\wide@bar@[3]{
	\begingroup
	\def\mathaccent##1##2{
		\let\mathaccent\save@mathaccent
		\if#32 \let\macc@nucleus\first@char \fi
		\setbox\z@\hbox{$\macc@style{\macc@nucleus}_{}$}
		\setbox\tw@\hbox{$\macc@style{\macc@nucleus}{}_{}$}
		\dimen@\wd\tw@ \advance\dimen@-\wd\z@ \divide\dimen@ 3 \@tempdima\wd\tw@ \advance\@tempdima-\scriptspace \divide\@tempdima 10 \advance\dimen@-\@tempdima \ifdim\dimen@>\z@ \dimen@0pt \fi \rel@kern{0.6}\kern-\dimen@
		\if#31 \overline{\rel@kern{-0.6}\kern\dimen@\macc@nucleus\rel@kern{0.4}\kern\dimen@} \advance\dimen@0.4\dimexpr\macc@kerna \let\final@kern#2 \ifdim\dimen@<\z@ \let\final@kern1 \fi
		\if \final@kern1 \kern-\dimen@ \fi
		\else \overline{\rel@kern{-0.6}\kern\dimen@#1} \fi }
	\macc@depth\@ne	\let\math@bgroup\@empty \let\math@egroup\macc@set@skewchar 	\mathsurround\z@ \frozen@everymath{\mathgroup\macc@group\relax} 	 \macc@set@skewchar\relax \let\mathaccentV\macc@nested@a	\if#31 \macc@nested@a\relax111{#1} \else \def\gobble@till@marker##1\endmarker{} \futurelet\first@char\gobble@till@marker#1\endmarker \ifcat\noexpand\first@char A\else \def\first@char{} \fi \macc@nested@a\relax111{\first@char} \fi
	\endgroup }
\makeatother

\newcommand{\PP}[1]{\mathbb{P}\left(#1\right)}

\newcommand{\EE}[1]{\mathbb{E}\!\left(#1\right)}

\newcommand{\bb}{\textbf{b}}

\newcommand{\bx}{\textbf{x}}
\newcommand{\by}{\textbf{y}}

\newcommand{\bX}{\textbf{X}}

\newcommand{\fp}{\mathfrak{p}}
\newcommand{\fP}{\mathfrak{P}}

\newcommand{\mM}{\mathcal{M}}
\newcommand{\mS}{\mathcal{S}}
\newcommand{\mH}{\mathcal{H}}
\newcommand{\mP}{\mathcal{P}}

\newcommand{\mG}{\mathcal{G}}

\newcommand{\mA}{\mathcal{A}}

\newcommand{\mQ}{\mathcal{Q}}

\newcommand{\RR}{\mathbb{R}}
\newcommand{\NN}{\mathbb{N}}
\newcommand{\ZZ}{\mathbb{Z}}
\newcommand{\II}{\mathbb{I}}

\newcommand{\rank}[1]{\text{rk} ( {#1} )}
\newcommand{\Sub}[2]{{#1} [{#2}]} 
\newcommand{\Vol}[1]{\mathrm{Vol} \left( {#1} \right)}

\definecolor{darkred}{cmyk}{.3,.9,.80,.2}

\title[Threshold functions for systems of equations in random sets]{Threshold functions and Poisson convergence for systems of equations in random sets}

\author[J. Ru\'{e}]{Juanjo Ru\'{e}}
\address{J. Ru\'{e}: Department of Mathematics, Universitat Polit\`ecnica de Catalunya and Barcelona Graduate School of Mathematics, Edifici Omega, 08034 Barcelona, Spain}
\email{juan.jose.rue@upc.edu}

\author[C. Spiegel]{Christoph Spiegel}
\address{C. Spiegel: Department of Mathematics, Universitat Polit\`ecnica de Catalunya and Barcelona Graduate School of Mathematics, Edifici Omega, 08034 Barcelona, Spain}
\email{christoph.spiegel@upc.edu}

\author[A. Zumalac\'arregui]{Ana Zumalac\'arregui}
\address{A. Zumalac\'arregui: Department of Pure Mathematics, University of New South Wales, 2052 NSW, Australia}
\email{a.zumalacarregui@unsw.edu.au}

\thanks{J.~R. was partially supported by the FP7-PEOPLE-2013-CIG project CountGraph (ref. 630749), the Spanish MICINN projects MTM2014-54745-P and MTM2014-56350-P, the DFG within the Research Training Group \emph{Methods for Discrete Structures} (ref. GRK1408), and the Berlin Mathematical School. C.~S. was supported by a Berlin Mathematical School Scholarship. A.~Z. is supported by the Australian Research Council Grant DP140100118. The first and third authors started this project while financed by the MTM2011-22851 grant (Spain) and the ICMAT Severo Ochoa Project SEV-2011-0087 (Spain).}

\begin{document}

\begin{abstract}

We study threshold functions for the existence of solutions to linear systems of equations in random sets and present a unified framework which includes arithmetic progressions, sum-free sets, $B_{h}[g]$-sets and Hilbert cubes. In particular, we show that there exists a threshold function for the property \textit{``$\mA$ contains a non-trivial solution of $M\cdot\bx=\textbf{0}$''} where $\mA$ is a random set and each of its elements is chosen independently with the same probability from the interval of integers $\{1,\dots,n\}$. Our study contains a formal definition of trivial solutions for any linear system, extending a previous definition by Ruzsa when dealing with a single equation.

Furthermore, we study the distribution of the number of non-trivial solutions at the threshold scale. We show that it converges to a Poisson distribution whose parameter depends on the volumes of certain convex polytopes arising from the linear system under study as well as the symmetry inherent in the structures, which we formally define and characterize.
\end{abstract}

\maketitle

\section{Introduction}\label{introduction}

The study of the existence or absence of a given configuration in a large combinatorial structure plays a central role in discrete mathematics and particularly in combinatorial number theory.

On one hand, the extremal aspects of this question have provided an active area of research in extremal combinatorics where many different techniques are exploited to obtain results. One major and well known example is Szemerédi’s Theorem~\cite{Sz75} which proves the existence of arbitrarily long arithmetic progressions in sets of positive density, see also~\cite{Gow2001,Fu77}. However, most of these results rely on smart ad hoc arguments that strongly depend on the specific structures under consideration. For some general results that have been obtained in the extremal case see Ruzsa~\cite{Rz93, Rz95} as well as Shapira~\cite{Sh06}.

On the other hand, the common behaviour, i.e. what to expect for most sets, has been less studied. In this scenario one is interested in obtaining results regarding the existence (or absence) of solutions to systems of linear equations in random sets. This paper aims to provide answers to this question, giving a clear picture for a wide variety of structures that have been studied in the extremal case.

The model for random sets we consider in this work is known as the \emph{binomial model} $[n]_p$ and is analogous to the $\mG(n,p)$  model in random graphs, which is defined with respect to the parameter $n$ and the probability $p$ (possibly depending on $n$).  We consider random sets $\mA$ where every element  $a\in\{1,\dots,n\}=[n]$ is chosen to be in $\mA$ independently with probability $\PP{a\in \mA}=p$.  In this context, we say that a certain property holds ``for almost every set $\mA$'' or that it holds ``asymptotically almost surely'' if the probability that a random set in $[n]_p$ does not satisfy such property tends to zero as $n$ tends to infinity.

Let us note that in $[n]_p$ any specific set $\mA\subseteq [n]$ will appear with probability
\[\mathbb{P}\,(\mA)= p^{|\mA|}(1-p)^{n-|\mA|},\]
thus all sets with the same size are equiprobable. Also, the expected size for a random set in this model is
\[\EE{|\mA|}=\sum_{k=1}^n \PP{a\in \mA} = np.\]
It follows that with high probability the number of elements in a random set in the described model will be close to $np$.

Let $P$ be a combinatorial property and $\mA$ a binomial random set in $[n]$ with parameter $p$. In this context, we say that $t(n)$ is a \emph{threshold} for the property $P$ if
\begin{itemize}
\item [(i)]  $p=o(t(n))$ implies  $\lim_{n\to\infty} \mathbb{P}\,(``\mA \text{ satisfies property } P\text{''})\rightarrow 0$ (the 0-statement) and
\item [(ii)]  $p=\omega(t(n))$ implies $\lim_{n\to\infty} \mathbb{P}\,(``\mA \text{ satisfies property } P\text{''}) \rightarrow 1$ (the 1-statement).
\end{itemize}
Roughly speaking: when $p$ is ``above the threshold'' almost all sets satisfy property $P$ and ``below the threshold'' almost no set satisfies the property. Observe that thresholds are only uniquely defined within constant factors. Therefore, we are interested in the order of magnitude of this transition phase.

\SP

Consider the homogeneous linear system of $r$ equations in $m$ variables (with $r<m$)
\begin{equation*}
\left\{\begin{array}{ccc}a_{11}x_1+\cdots+a_{1m}x_m & = & 0 \\
& \vdots & \\
a_{r1}x_1+\cdots+a_{rm}x_m & = & 0
\end{array}{}\right.
\end{equation*}
which will be identified with its corresponding matrix $M = (a_{ij})\in \mM_{r,m}(\ZZ)$. We will only consider matrices $M = (a_{ij})\in \mM_{r,m}(\ZZ)$ satisfying three natural conditions:
\begin{enumerate}
\item \emph{positivity} (the system must contain at least one solution with all positive entries),
\item \emph{irredundancy} (for each $i \neq j$ there exists a solution $\bx = \left(x_1,\dots,x_m\right)$ such that $x_i \neq x_j$),
\item \emph{non-degeneracy} (the matrix has full rank).
\end{enumerate}
From now on we call any system of equations given by a matrix satisfying these conditions \emph{admissible}. Roughly speaking, an admissible system of equations must have positive solutions without repeated coordinates and cannot be reduced to another one with a smaller number of equations or variables.

For a given set $\mA\subseteq [n]$, we say that $\mA$ \emph{contains a solution to} $M\cdot \bx=0$ if there exist elements $x_1,\cdots, x_m\in \mA$ such that $\bx = (x_{1},\cdots,x_{m})^\textrm{T}$ is a solution to the system.

\SP

The problem we address in this paper is the following one: let $\mA \subseteq [n]$ be a random set sampled from $[n]_p$. For an admissible matrix $M$ we study how the random variable
$$\bX = |\mA^m\cap \{\bx \text{ non-trivial} \,:\,  M\cdot \bx=\textbf{0}\}|$$
behaves with respect to $p$ and deduce the existence of a threshold function for the combinatorial property \textit{``$\mA$ contains a non-trivial solution of $M \cdot \bx=\textbf{0}$''}. The existence of such a function is assured by the fact that monotone properties of random sets always have thresholds functions~\cite{BT87}. Note that a crucial point in the analysis of the solutions is the definition of \textit{trivial solutions}. It might be clear what a trivial solution looks like for $k$-AP or a $B_h[g]$-set, but in the general setting this concept is less obvious and will be explored in Section~\ref{trivialsol}. This will be the main difference compared to the study of small subgraphs in the classical $\mG(n,p)$ model.

\SP

The exponent of the threshold function will have to be maximized over all \emph{induced submatrices} of the system $M \cdot \bx=\textbf{0}$. A full explanation will follow in Section~\ref{submat}, but in order to state our results we need to introduce a definition here. For any set of column indices $\emptyset \subseteq Q \subseteq [m]$ let $M^Q$ denote the matrix obtained from $M$ by keeping only the columns indexed by $Q$, where $M^{\emptyset}$ is the empty matrix. We write the rank of a matrix $M$ as $\rank{M}$ and define $r_Q = r - \rank{M^{\widebar{Q}}}$ for all $\emptyset \neq Q \subseteq [m]$. Here we let $\rank{M^{\emptyset}} = 0$. This allows us to define the following parameter.

\begin{definition}\label{def:maxparameter}
For any admissible matrix $M\in \mM_{r,m}(\ZZ)$ define
\begin{equation}
	c(M) = \max_{\emptyset\neq Q\subseteq\left[m\right]} \frac{|Q|}{|Q|-r_Q}.	
\end{equation}
\end{definition}

Maximizing a parameter over all possible induced submatrices is reminiscent of similar results obtained while studying analogous problems in graph theory~\cite{ER60}. The intuition behind it will be explored later with some examples. Also note the similarity between this parameter and the one developed by R\"odl and Ruci\'nski, who restricted themselves to density regular matrices when looking at random sparse versions of Rado's theorem in~\cite[Definition~1.1]{RR97}. Under the previous assumptions and using these definitions, we obtain the following result:

\begin{theorem}\label{thm:threshold}
For some $r<m$, let $M\in \mM_{r,m}(\ZZ)$ define an admissible system. Then, the probability $p=n^{-1/c(M)}$ is a threshold function for the property: ``$\mA$ contains a non-trivial solution of $M \cdot \bx=\textbf{0}$''.
\end{theorem}

In other words, whenever the size of $\mA$ is $o\big(n^{1 - 1/c(M)}\big)$ we can assure that asymptotically almost surely there are no solutions, other than trivial ones, of the linear system $M\cdot \bx=\textbf{0}$ with $\bx \in \mA^{m}$. The main contribution in the study will come from \emph{proper} solutions, i.e. those solutions with pairwise different coordinates, since roughly speaking solutions with repeated coordinates only start appearing for larger values of $p$.

We also study the behaviour of the limiting probability at the threshold. With this purpose in mind, observe that any system of equations $M\cdot \bx=0$ together with the restrictions $x_i\in[0,1]$ define a non-empty, convex and rational polytope of dimension $m-r$ which we denote by $\mP_M$. We show that if and only if our system is strictly balanced (see the corresponding definition in Section~\ref{submat}) and $p=Cn^{-1/c(M)}$ for some constant $C>0$, the limiting distribution converges to a Poisson distribution and hence the probability of having a solution tends to $1$ with an exponential decay in $C$. Furthermore, the parameter only depends on the volume of $\mP_{M}$ as well as the inherent symmetry of the system (see Definition~\ref{def:symmetryconstant}). More precisely,
\begin{theorem}\label{thm:local}
For some $r<m$, let $M \in \mM_{r,m}(\ZZ)$ define an admissible system, $p=C n^{-1/c(M)}$ for some $C>0$, $\mu= \frac{\Vol{\mP_{M}}}{\sigma(M)} C^m $. Then for every non-negative integer $t$
\begin{equation*}
	\lim_{n\rightarrow \infty}\PP{``\mA \text{ contains } t \text{ non-trivial solutions of } M \cdot \bx=\textbf{0}\text{''}}=\tfrac{\mu^t}{t!}e^{-\mu},
\end{equation*}
if and only if  $M$ is strictly balanced. Here $\mP_{M}$ is the polytope associated to the system $M\cdot \bx=\textbf{0}$ and $\sigma(M)$ is a computable constant depending on $M$.
\end{theorem}
Observe that the previous result implies that for $p=Cn^{-1/c(M)}$ the number of solutions is approximately Poisson distributed with parameter $\mu$, and, in particular
$$\lim_{n\rightarrow \infty}\PP{``\mA\text{ contains a non-trivial solution of }M\cdot \bx=\textbf{0}"}=1-e^{-{\mu}}.$$

Note also that for the statement of Theorem~\ref{thm:threshold} it was not of importance whether one regards the solutions as subsets of $\mA$ or as vectors in $\mA^m$. However, for the behaviour at the threshold, one needs to be more careful since the constants that appear when counting solutions play a role. This is not an issue when dealing with the usual systems considered in the literature, such as $k-$AP, Sidon and $B_h[g]$-sets (see below). However, since we are developing our theory for \emph{any} admissible system, we have to take greater care of this issue since one can construct specific examples where neither approach aligns with one's intuition. Our approach will therefore be somewhere in between considering solutions as vectors or as subsets. We start out considering them as vectors and then take care of a symmetry factor which depends only on $M$ and that can be unilaterally applied to all solutions. Note that for the already mentioned common systems there is no difference between this approach and considering them as subsets. More details are given in Section~\ref{sec:countingsolutions}, including examples that illustrate why this approach is prudent.

The computation of the constant $\Vol{\mP_M}$ appearing in Theorem~\ref{thm:local} is an algorithmically involved  problem when dealing with general systems. One could compute this volume by means of triangulations of the polytope~\cite{LRS10}, but the problem is in general (for dimension greater than~$3$) NP-complete~\cite{BEF00}. We provide computations for some concrete systems in Section~\ref{c(M)}.

\SP

This work will include a precise analysis of interesting combinatorial families that have been studied in the literature from many different points of view and which fit into the presented scheme. Let us state some of these common configurations:

\subsubsection*{Arithmetic progressions} A set of integers is an arithmetic progression of length $k$ (or shortly, a $k$-AP) if it can be written in the form $a,\,a+d,\,\dots,a+(k-1)d$ for some $a,\, d \in \mathbb{Z}$ and $d \neq 0$.

\subsubsection*{Sidon and $B_h[g]$-sets} A set of integers $\mA$ is called a \emph{Sidon set} (or $B_2[1]$-set) if every integer $k$ has at most one representation as a sum of two elements of $\mA$, up to permutations of the summands involved. One can generalize this concept in several ways; for example a set $\mA$ of non-negative integers is a $B_{h}[g]$-set if every integer has at most $g$ representations as a sum of $h$ elements of $\mA$, modulo permutations of the summands involved.

\subsubsection*{Hilbert Cubes} Another possible generalization of Sidon sets are so-called Hilbert cubes: a set $\mH$ of integers is a Hilbert cube of dimension $k$ (or $k$-cube) if there exist positive and distinct integers $h_0,h_1,\dots,h_k$ satisfying
$$\mH = \Big\{h_0+\sum_{i=1}^k \epsilon_i h_i:\, \epsilon_i\in\{0,1\}\Big\}.$$
Clearly a set $\mA$ is Sidon if it does not contain any $2$-cube. As it is shown by S\'andor~\cite{Sa07} almost all sets in $[n]$ with size $\omega\big(n^{1-(k+1)/2^k}\big)$ contain a $k$-cube. Our results extend those in~\cite{Sa07}, proving in particular Conjecture 3.1 in the same work.

\subsubsection*{Barycentric and sum-free sets} A set $\mA$ contains a $k$-barycentric set if there exist elements $a_1,\dots,a_k,a_{k+1}\in \mA$ such that $$a_1+a_2+\dots+a_k=ka_{k+1},$$ that is $a_{k+1}$ is the average of $a_1,\dots,a_k$. Clearly if $k=2$ we have a $3$-AP and trivial solutions are given by $a_1=\cdots=a_{k+1}$. Finally, a set of integers $\mA$ is a $k$-sum-free set if for every pair $a,\,a'\in \mA$ the sum $a+a'$ is not an element of $k \cdot \mA = \{ka: a\in\mA\}$. The case $k=1$ is also known as \emph{sum-free} equation or \emph{Schur} equation.

\SP

The existence of such structures can be expressed using systems of equations of the type $M\cdot \bx=0$ for admissible matrices $M\in \mM_{r,m}(\ZZ)$. A set $\mA$ avoids a $k$-AP if the homogeneous system defined by
\[M_{k\text{-AP}}=\left(
\begin{smallmatrix}
1&-2&1&  \\
      &1&-2&1&&& \\
      &&&&\cdots\\
      &&&&&1&-2&1
\end{smallmatrix}
\right)\in \mM_{k-2,k}(\ZZ),\]
does not have a non-trivial solution $\bx=(x_1,\dots, x_{k})\in \mA^{k}$. In this case all trivial solutions (see Definition~\ref{def:trivial}) are given by $x_1=x_2=\cdots=x_k\in \mA$ and correspond to the case $d=0$.
A set $\mA$ is a Sidon set if there are no solutions $\bx=(x_1,x_2,x_3,x_4)\in \mA^4$ of the linear system $x_1+x_2=x_3+x_4$, except for the trivial ones, which have the form either $(a,b,a,b)$, $(a,b,b,a)$ for $a,b\in [n]$. Similarly, a set $\mA$ is $B_h[g]$ if there are no solutions in $\mA^{h(g+1)}$ of the linear system defined by
\begin{equation}\label{eq: matrix_Bhg}
M_{B_h[g]} =\left(\begin{smallmatrix}
      1&\stackrel{h}{\cdots}&1&-1&\stackrel{h}{\cdots}&-1 \\
      &&&1&\stackrel{h}{\cdots}& 1&-1& \stackrel{h}{\cdots}& -1\\
      &&&&&&&&&\cdots\\
      &&&&&&&&&&1&\stackrel{h}{\cdots}&1&-1&\stackrel{h}{\cdots}&-1
\end{smallmatrix}\right)\in \mM_{g,h(g+1)}(\ZZ).
\end{equation}
A set $\mA$ avoids $3$-Hilbert cubes if it does not contain solutions to
$$M_{\mH_3}=\left(\begin{smallmatrix}
-1&1&1&-1  \\
&&-1&1&1&-1\\
&&&&-1&1&1&-1\\
&-1&1&&&-1&1 \\
\end{smallmatrix}\right)\in \mM_{4,8}(\ZZ),$$
and in general for a $k$-Hilbert cube we will have a system of rank $2^k-(k+1)$ in $2^k$ variables. Lastly, a set $\mA$ is $k$-sum-free if there are no solutions $\bx=(x_1,x_2,x_3)\in \mA^3$ of the linear system $x_1+x_2=kx_3$, (when $k=2$ the problem corresponds with a 3-AP).

It is clear from the definition of $M_{k\text{-AP}}$, $M_{B_h[g]}$ and the $k$-sum-free family that all matrices have maximum rank, that is $r=k-2$, $r=g$ and $r=1$, respectively and are in fact admissible matrices. The application of the previous theorems and the computations in Section~\ref{c(M)} are summarized in Table~\ref{Table}.

\begin{table}[htb]\label{Table}
\begin{center}
\begin{tabular}{c|cc|cc|cc}
	& $r$ & $m$ & $p$ & $\EE{|\mA|}$ & $\Vol{\mP_M}$&$\sigma(M)$ \Bstrut \\
	\hline
	$k\mathrm{-AP}$ & $k-2$  & $k$  & $n^{-2/k}$  &  $n^{1-2/k}$  & $1/(k-1)$ & $2$ \Tstrut \\
	$\mathrm{Sidon}$ & $1$  & $4$ & $n^{-3/4}$  & $n^{1/4}$ & $2/3$ & $8$ \Tstrut  \\
	$B_{h}[g]$ & $g$ & $h(g+1)$  & $n^{\frac{g}{h(g+1)}-1}$ & $n^{\frac{g}{h(g+1)}}$ & Section~\ref{c(M)}&$(g+1)!(h!)^{g+1}$ \Tstrut  \\
	$k\mathrm{-cube}$ & $2^k-(k+1)$ & $2^k$ & $n^{-\frac{k+1}{2^k}}$ & $n^{1-\frac{k+1}{2^k}}$ & $\frac{2^{2^k-1}}{(k+1)!k!}$ & $2^{2^k-1}$ \Tstrut  \\
	$\mathrm{sum-free}$ & $1$ & $3$ & $n^{-2/3}$  & $n^{1/3}$ & $1/2$ & $2$ \Tstrut  \\
	$k\mathrm{-sum-free}$ & $1$ & $3$ & $n^{-2/3}$  & $n^{1/3}$ & $1/k$ & $2$ \Tstrut  \\
	$k\mathrm{-barycentric}$ & $1$ & $k+1$ & $n^{-k/k+1}$ & $n^{1/k+1}$ & $1/k$ & $k!$ \Tstrut  \\
\end{tabular}
\bigskip
\caption{Threshold for different combinatorial families. The analysis of $k$-cubes is fully done in \cite{Sa07}.}
\label{table:thresholds}
\end{center}
\end{table}
Nevertheless, let us remark that the general approach presented in this paper allows one to study a lot more than just these few linear structures. Therefore, some generalizations had to be made that proved to be very delicate when no longer considering these well-known examples. We have already mentioned the issue of defining trivial solutions as well as the symmetry between solutions. We have included the computation of the \emph{symmetry constant} $\sigma(M)$ in the previous table (see Section~\ref{sec:countingsolutions} for a proper definition). We will also need to count the number of solutions to the system in $[n]^m$ and, as we will see in Sections~\ref{sec:countingsolutions} and~\ref{submat}, one must be very careful when doing so for a general system.

\SP
\subsection*{Outline.} In \emph{Section~\ref{sec:extremalcase}} we compare the results obtained from Theorem~\ref{thm:threshold} to what is currently known in the extremal case for certain structures. \emph{Section~\ref{tools}} contains a brief overview of the tools needed for the proofs later on in this paper. In \emph{Section~\ref{sec:countingsolutions}} we will apply Ehrhart's Theory to obtain a useful lemma and a simple corollary that allow us to count the number of solutions to a system of linear equations up to some symmetry. \emph{Section~\ref{trivialsol}} contains a formal definition of trivial solutions with examples to motivate it. In \emph{Section~\ref{submat}} we introduce the concept of induced submatrices and an important proposition regarding non-trivial solutions. It also contains a formal definition of strictly balanced systems which is a necessary prerequisite for Theorem~\ref{thm:local}. A proof of Theorem~\ref{thm:threshold} using the Second Moment Method can be found in \emph{Section~\ref{sec:thresholdproof}} and in \emph{Section~\ref{Brun-proof}} we study the local behavior of the threshold which results in a proof of Theorem~\ref{thm:local} via an application of Brun's Sieve. The analysis of $\Vol{\mP_{M}}$ associated to certain combinatorial families is carried out in \emph{Section~\ref{c(M)}}. Finally, in \emph{Section~\ref{related}} we discuss related problems and generalizations.

\section{State of the art in the extremal case} \label{sec:extremalcase}

In the presented approach, we seek to give a picture of the qualitative behavior of a random set. However one might wonder how far the common situation is from the extremal cases. The problem of estimating the size of maximal sets avoiding the structures introduced above has been intensively studied. In this direction one can find several results which give upper bounds for sets avoiding specific structures or, on the opposite direction, explicit constructions of large sets with this property. In both cases one requires \emph{ad hoc} arguments that strongly depend on each specific structure.

\subsubsection*{Arithmetic progressions}

For sets avoiding $k$-AP's we must go back to Szemer\'edi's Theorem, that states that no set with positive density can avoid $k$-AP's for any $k$. In particular, for $k=3$ non-trivial bounds were first obtained by Roth~\cite{Ro53} and then refined by several authors, see~\cite{HeBr87,Bou08}. Nowadays, the best upper bound is established by Bloom~\cite{Bl16}. On the other hand, Behrend~\cite{Be46} constructed a set avoiding $3$-AP's of large size; this construction was slightly improved by Elkin~\cite{Elk11} (see also~\cite{GW10}). More precisely, we have
$$ n\cdot \frac{(\log n)^{1/4}}{e^{c\sqrt{\log n}}} \ll\max_{\mA\subset [n]}\{|\mA|:\, \mA \text{ avoids }3\text{-AP's }\}\ll  n\cdot \frac{(\log\log n)^4}{\log n}, $$
for some constant $c$.

Concerning the general $k$-AP problem, analogous bounds have been obtained: the upper bounds come from the pioneering work of Gowers~\cite{Gow2001} and, more recently, dense constructions that lead to lower bounds for this problem were stablished by O'Bryant~\cite{OBr11}. These results can be summarized as follows
$$n\cdot \frac{(\log n)^{(2\log k)^{-1}}}{e^{c(k)(\log n)^{\log^{-1} k}}} \ll\max_{\mA\subset [n]}\{|\mA|:\, \mA \text{ avoids }k\text{-AP's }\}\ll  n\cdot (\log\log n)^{-2^{-2^{(k+9)}}},$$
for a certain constant $c(k)$ only depending on $k$.

We show that almost all sets with size $\omega(n^{1-2/k})$ contain $k$-AP's. Observe that, for $k=3$ the gap between the usual situation and the extremal set is very large: most sets with size $\omega(n^{1/3})$ contain $3$-AP's but there are examples of (almost) linear size avoiding this structure. Nevertheless, as $k$ grows to infinity, this quantity approximates to $n$ and the gap between the exponents tends to $0$.

In the direction of the present article, Warnke has also studied the upper tail of the number of $k$-arithmetic progressions and Schur triples in random subsets, establishing exponential bounds~\cite{Wa16}.

\subsubsection*{Sidon and $B_h[g]$-sets}

The study of Sidon sets dates back to Erd\H{o}s. In~\cite{ET41} Erd\H{o}s and Tur{\'a}n obtained an upper bound for the size of a maximal Sidon set in $[n]$ (see~\cite{Li69, Cil10} for further improvements of this result). In fact, there are algebraic constructions of Sidon sets that, combined with Erd\H{o}s-Tur{\'a}n result, prove
$$\max_{\mA\subset [n]}\{|\mA|:\, \mA \text{ is Sidon }\}\sim n^{1/2}.$$
In the direction of the present article, the $B_h[1]$ case was studied in detail by Godbole, Janson, Loncatore and Rapoport in~\cite{GJLR99}. They show that almost no set with size $\omega(n^{1/2h})$ is $B_h[1]$. The proof is based on a tailor-made analysis on the particular shape of the equations defining $B_h[1]$ sets.

Clearly, for $h=2$ (that is Sidon), the gap between the exponents in the usual situation, namely $|\mA|=o(n^{1/4})$, and the extremal one, say $|\mA|=n^{1/2}(1+o(1))$, is very big. Let us also mention that in~\cite{KLRS15} Kohayakawa, Lee, R\"odl and Samotij study the number of Sidon sets and the maximum size of Sidon sets contained in a sparse random set of integers. In particular, in Section 5 they analyze, by means of  the Kim-Vu polynomial concentration inequality~\cite{KV00}, the number of solutions to the Sidon equation (when the probability lies above the threshold). These results could be deduced using the presented framework. Concerning the general case, it is known that the cardinality of a maximum $B_h[g]$-set in $[n]$ is $\asymp n^{1/h}$, but the main difficulty is to obtain a precise constant for the problem~\cite{CRT02, CRV10}. As we show in Theorem~\ref{thm:threshold} almost all sets in $[n]$ of size $o(n^{{1}/{h}-{1}/{h(g+1)}})$ are $B_{h}[g]$. Once again, if we fix $h$ and let $g$ grow to infinity both situations approach each other.

\subsubsection*{Hilbert Cubes}

Hilbert originally proved that any finite coloring of the positive integers contains a monochromatic $k$-cube. The density version of this result is known as Szemer\'edi's Cube Lemma and it is a key point in his proof of Roth's Theorem. Gunderson and R\"{o}dl~\cite{GR98} obtained, by counting arguments, that for sufficiently large $n$, any set $\mA\in [n]$ with size $2n^{1-1/2^{k-1}}$ contains a $k$-cube. On the other side, by means of probabilistic arguments, one can construct a set of size $n^{1-k/(2^{k}-1)}$ avoiding $k$-cubes. For the particular case $k=3$, very recently Cilleruelo and Tesoro~\cite{CT15} have obtained an algebraic construction of a set of size $\gg n^{2/3}$ avoiding $3$-cubes.

As in the previous cases, when $k$ grows the existing gap between the exponents in our result and the ones in the upper and lower bounds tends to $0$.

\SP

\subsubsection*{$k$-sum free sets}

The question of maximizing the cardinality of a set of integers in $[n]$ avoiding $x+y=kz$ belongs to the folklore: one cannot select more than $\lceil\frac{n}{2} \rceil$ integers satisfying this condition for $k=1$ and this is optimal. The case $k=2$ coincides with the exclusion of $3$-AP's. Concerning $k=3$, the problem was solved by Chung and Goldwasser~\cite{CG96_2} getting the same estimates as for $k=1$. For $k\geq 4$, and sufficiently large $n$, Chung and Goldwasser~\cite{CG96_1} discovered $k$-sum-free sets of linear size in $n$ (and density tending to $1$ as $k$ increases); in fact Baltz, Hegarty, Knape, Larsson and Schoen~\cite{BHKLS05} showed that this construction is optimal. Therefore, for this family it is known that the maximal size of a $k$-sum-free set is linear in $n$ but Theorem~\ref{thm:threshold} asserts that almost all sets of size $\omega(n^{1/3})$ contain at least one solution to $x+y=kz$, for every $k$. Observe that in this family, the parameter $k$ does not play a role in the position of the  threshold.

\section{Tools}\label{tools}

In this section we recall the Second Moment Method, Janson's inequality and Brun's Sieve – in the context of the Probabilistic Method – as well as basic notions in Ehrhart's Theory for counting lattice points in convex polytopes.

\subsection{The Second Moment Method.}\label{sub:2nd-moment}
The Second Moment Method is used in the version given by Corollary 4.3.4. of Alon and Spencer \cite{AS08}: let $X=\II_1+\dots+\II_s$ be a sum of $s = s(n)$ indicator random variables, where $\II_i = \II_i(n)$ corresponds to some event $E_i = E_i(n)$. For convenience we suppress the dependence on $n$ which defines our $o, O$ notation. We write $i \sim j$ if $i\neq j$ and the events $E_i$ and $E_j$ are not independent. Define
\begin{equation}\label{eq:delta}
	\Delta = \sum_{i\sim j} \PP{E_i \wedge E_j}
\end{equation}
If $\EE{X} \rightarrow \infty$  and $\Delta=o\big( \, \EE{X}^2\big)$ (as $s\rightarrow \infty$), then $X\sim \EE{X}$ asymptotically almost surely. In particular, under these assumptions, $X>0$ with probability tending to $1$.\\

\subsection{Brun's Sieve}\label{Brun}
The traditional approach to the Poisson Paradigm is used in the version given by Theorem 8.3.1. of Alon, Spencer~\cite{AS08}: let $X=\II_1+\dots+\II_s$ again be a sum of $s$~indicator random variables associated to some events $E_i$. Let
$$S^{(t)}=\sum_{\{i_1,\dots, i_t\} \in \binom{[s]}{t}} \PP{E_{i_1}\wedge E_{1_2}\wedge\cdots\wedge E_{i_t}},$$
where the sum is taken over all subsets $\{i_1,\dots, i_t\}\subseteq [s]$ of $t$ elements. The Inclusion-Exclusion Principle gives us
$$\PP{X=0} = \PP{\widebar{E_1}\wedge \cdots \wedge \widebar{E_s}}=1-S^{(1)}+S^{(2)}-\cdots + (-1)^tS^{(t)}+\cdots$$

\begin{theorem}[Brun's Sieve]\label{Thm:Brun}
Suppose there is a constant $\mu$ so that $\EE{X}=S^{(1)}= \mu(1+o(1))$, and such that for every fixed $t$
$$S^{(t)}\longrightarrow \frac{\mu^t}{t!}.$$
Then, for every non-negative integer $t$
$$\PP{X=t}\longrightarrow \frac{\mu^t}{t!}e^{-\mu}.$$
\end{theorem}

\subsection{Lattice points in dilates of polytopes – Ehrhart's Theory.}\label{sub:ehrart-theory} 

A basic reference for definitions and first properties of convex polytopes is~\cite{Zi95}. For further results in lattice points in rational polytopes, see~\cite{BR07,Lo05}.

A \emph{convex polytope} is the convex hull of a finite set of points (which are always bounded), or a bounded intersection of a finite set of half-spaces, and is said to be \emph{rational} (resp. \emph{integral}) if its vertices (i.e. its corner points) are points with rational (resp. integral) coordinates. Every rational polytope has a matrix representation of the form
\begin{equation} \label{eq:polytopetomatrix}
\{\textbf{x}\in \mathbb{R}^k: P \cdot \textbf{x} \geq \textbf{b}\},\, P \in \mM_{l,k}(\ZZ),\,\textbf{b}\in \ZZ^k
\end{equation}
for certain non-negative integers $m,d$. Note that the inequalities can be easily turned into equalities through the use of slack variables. The \emph{(relative) dimension} of a polytope $\mP$ is the dimension of the affine space $\text{span} \, \mP := \left\{ \bx + \lambda (\by - \bx) \: : \: \bx, \by \in \mP , \: \lambda \in \RR \right\}$. Note that this dimension is not necessarily $d$, but a smaller non-negative integer.  For a given polytope $\mP$, let $\Vol{\mP}$ be the volume of $\mP$ in this affine space and $n\cdot \mP=\{n\mathbf{p}: \mathbf{p}\in \mP\}$ the $n$th-dilate of the polytope.

Ehrhart's Theorem~\cite{Eh62} (see also~\cite{Mac63}) gives a precise description of the number of integer points on the $n$th-dilate of a rational polytpe in this context: the quantity $ \left|n\cdot \mP \cap\, \mathbb{Z}^{\mathrm{dim}(\mP)}\right|$ is given by a pseudopolynomial in $n$ of degree $\mathrm{dim}(\mP)$ (recall that a pseudopolynomial is a function $f(n)=\sum_{i=0}^d c_i(n)n^i$ where the functions $c_0(n),\dots, c_d(n)$ are periodic). More precisely, we have the following theorem:
\begin{theorem}[Ehrhart's Theorem]\label{thm:ethm}
 Let $\mP$ be a $d$-dimensional convex polytope.
\begin{itemize}
\item[i.] If $\mP$ is an integral polytope, then $\left|n\cdot \mP \cap\, \mathbb{Z}^{d}\right|$ is a \emph{polynomial} in $n$ of degree $d$.

\item [ii.] If $\mP$ is a rational polytope, then $\left|n\cdot \mP \cap\, \mathbb{Z}^{d}\right|$ is a \emph{pseudopolynomial} in $n$ of degree $d$. Additionally, its period divides the least common multiple of the denominators of the coordinates of the vertices of $\mP$.
\end{itemize}
\end{theorem}

Additionally to Theorem~\ref{thm:ethm}, one can easily show that the leading coefficient in both cases is equal to $\Vol{\mP}$. As a trivial corollary, for a rational polytope $\mP$ of dimension $\mathrm{dim}(\mP)$ embedded in $\mathbb{R}^{\mathrm{dim}(\mP)}$, we have
\begin{equation}\label{eq:leadingcoefficient}
 \left|n\cdot \mP \cap\, \mathbb{Z}^{\mathrm{dim}(\mP)}\right|= \Vol{\mP} n^{\mathrm{dim}(\mP)}(1+o(1)).
\end{equation}
Let us mention that the full version of Ehrhart's Theorem will be used in Subsection~\ref{ssec:Bh} in order to study volumes in $B_h[g]$ families. However, the weaker version stated in Equation \eqref{eq:leadingcoefficient} will be use in the rest of the proofs.

\section{Counting proper solutions up to symmetry}\label{sec:countingsolutions}

Consider the system defined by $M \in \mM_{r,m}(\ZZ)$ and recall that a solution to it is called proper if its coordinates are pairwise different. Note that the number of proper solutions in $[n]$ to the given system has the obvious upper bound $n^{m-r}$. In fact, this bound trivially holds even for inhomogeneous systems, that is for any $\bb \in \ZZ^r$ and $n \in \NN$ we have
\begin{equation}\label{eq:trivialupperbound}
	\big| \{ \bx \in [n]^m : M \cdot \bx = \bb \} \big|	\leq n^{m-r}.
\end{equation}
One can also easily give an easy constructive proof for the existence of some constant $C > 0$ such that $Cn^{m-r}$ is a lower bound for $n$ large enough, see for example~\cite{JR11}. These two bounds are in fact sufficient to prove the statements in Theorem~\ref{thm:threshold}.

However, the statement in Theorem~\ref{thm:local} requires the exact asymptotic value of the ratio of the number of solutions to $n^{m-r}$. This is where we apply Ehrhart's Theory and in particular Equation~\eqref{eq:leadingcoefficient} in order to count the number of proper solutions to the system and obtain the fundamental Lemma~\ref{fundamental}. Solutions that have repeated coordinates will simply be considered as as proper solution to some reduced system as we will introduce in the next section.

Let us start by considering that two proper solutions which are counted as separate by Ehrhart's Theory can be essentially the same when considering symmetry. As an easy example for this consider that $3$-APs are given by $x_1 + x_3 = 2x_2$ for which $(1,2,3)$ and $(3,2,1)$ both are proper solutions that one might consider as essentially identical. However the situation is not quite as simple as grouping solutions together if they are identical up to permutation. Consider for example the system given by
\begin{equation}\label{eq:referee_example}
	x_1 + x_2 + x_3 = x_4 + x_5 + x_6 + x_7
\end{equation}
for which both $(2,3,100;1,4,49,51)$ and $(1,4,100;2,3,49,51)$ are again proper solutions. In this case however, one should consider them to be essentially different because the permutation did not just occur between coordinates with identical coefficients and hence cannot be applied to all solutions. The semicolons in the vector representations delineating the coordinates with factor $-1$ from those with factor $1$ were added to emphasize that fact. In order to deal with this distinction we denote by $\mS_m$ the symmetric group on $m$ elements and introduce the following definition:

\begin{definition}\label{def:symmetryconstant}
For any matrix $M \in \mM_{r,m}(\ZZ)$  its symmetry constant is defined as
$$\sigma(M) = | \left\{ \pi \in \mS_m \mid \text{all } \bx \text{ s.t. } M \cdot \bx = \mathbf{0} \text{ satisfy } M \cdot \bx_\pi = \mathbf{0} \right\}|,$$
where $\bx_\pi$ is the vector obtained by permuting the coordinates of $\bx$ according to $\pi$.
\end{definition}

This definition based on the solution space might be immediately applicable for systems such as the ones given in the introduction, where the solutions are intuitively clear. However if one is given a more complex and less structured system, the following simple characterization of the symmetry constant based solely on the matrix $M$ might be easier to apply:

\begin{lemma}\label{symmetrycharacterization}
	For any matrix $M \in \mM_{r,m}(\ZZ)$ and permutation $\pi \in \mS_m$ let $M_\pi$ denote the matrix obtained by permutating the columns of $M$ according to $\pi$. We have
	$$\sigma(M) = |\left\{ \pi \in \mS_m \mid M_\pi \cong M \right\}|,$$
	where $\cong$ denotes equality up to linear transformations of the rows.
\end{lemma}

\begin{proof}
The inequality $\sigma(M) \geq \left\{ \pi \in \mS_m  \mid M_\pi \cong M \right\}$ is trivial. In order to show equality, note that for every permutation $\pi$ from the set
$$\left\{ \pi \in \mS_m  \mid \text{all } \bx \text{ s.t. } M \cdot \bx = \mathbf{0} \text{ satisfy } M \cdot \bx_\pi = \mathbf{0} \right\}$$
we have $\text{ker}(M) \subseteq \text{ker}(M_{\pi^{-1}})$ and since both kernels have dimension $m-r$ we have equality. Now the kernel of a matrix is the orthogonal of the span of its rows and therefore the rows of $M_{\pi^{-1}}$ can be obtained by linear transformation of the rows in $M$.
\end{proof}

In the following we will consider two solutions $\bx,\by$ to be \textit{essentially different} if and only if
\[ \by\ne \bx_{\pi}\quad \text{for every}\quad \pi \in \mS_m \text{ satisfying } M_\pi \cong M, \]
and otherwise $\bx$, $\by$ are treated as the same solution, even if they might differ as vectors. Observe that the definition and its characterization captures the previous examples as intended, i.e. it considers $(1,2,3)$ and $(3,2,1)$ to be the same $3$-AP while distinguishing between two solutions given for the second example. For the common systems mentioned in the introduction this approach of considering solutions as vectors in $\mA^m$ and dividing by $\sigma(M)$ results in counting subsets in $\mA$ that are solutions. For the example~\eqref{eq:referee_example}, however there is a small but significant difference between this and our approach.

\SP

Now in order to apply Equation~\eqref{eq:leadingcoefficient} note that for any admissible matrix $M\in \mM_{r,m}(\mathbb Z)$ the system $M\cdot \textbf{x}=\mathbf{0}$ with the additional restraints $0\leq x_1,\dots, x_m\leq 1$ defines a rational polytope $\mP_M$ of dimension $m-r$ (by assumption the system has the maximum possible rank and the polytope is not empty by the positivity assumption). It is just the intersection of the $(m-r)$-dimensional solution space and the $m$-dimensional unit hypercube. Using this we can now formulate the following lemma which will be applied in the forthcoming sections and simplify the discussion:

\begin{lemma}\label{fundamental}
Let $r<m$, $M\in \mM_{r,m}(\mathbb Z)$ an admissible matrix and $\mP_M$ the rational polytope defined by $M$. Then the number of different proper solutions $\textbf{x}\in [n]^m$ of $M\cdot \textbf{x}=\mathbf{0}$ is of the form $\Vol{\mP_M}/\sigma(M) \: n^{m-r} \: (1+o(1))$.
\end{lemma}

\begin{proof}
The number of lattice points in $n\cdot \mP_M$ is precisely the number of (not necessarily proper) vector solutions to $M\cdot \textbf{x}=\mathbf{0}$ with the added condition that $\textbf{x}\in [n]^m$. As the intersection of the $(m-r)$-dimensional solution space and the $m$-dimensional unit hypercube, the polytope $\mP_M$ also has dimension $m-r$. By Equation~\eqref{eq:leadingcoefficient} the number of lattice points in the dilate $n\cdot \mP_M$ is simply $\Vol{\mP_M} n^{m-r}(1+o(1))$. Noting that we do not distinguish between solutions that are identical up to certain permutations as specified in Definition~\ref{def:symmetryconstant} introduces the factor $1/\sigma(M)$.

We therefore have to consider the set of solutions $\textbf{x}\in [n]^m$ of $M\cdot \textbf{x}=0$ with some repeated coordinates and show that they have a negligible contribution to the total number of solutions. These solutions belong to the intersection of $\mP_M$ with a subspace defined by repetitions of coordinates. Since by assumption our system is irredundant, $\mP_M$  contains at least one solution with no repeated coordinates; this implies that there is no subspace defined by the the repetition of coordinates containing $\mP_M$.  Therefore, the polytope resulting from the intersection has dimension strictly smaller than $m-r$.

It follows again by Equation~\eqref{eq:leadingcoefficient} that the number of solutions with certain repeated coordinates is $O(n^{m-r-1})$. Finally, the number of possible constellations of repeated coordinates is bounded by the number of partitions of $\{1,\dots,m\}$, so the total number of solutions with repeated coordinates is  $o(n^{m-r})$ and the lemma follows.
\end{proof}

\section{Trivial solutions}\label{trivialsol}

The key point of this section is to correctly define what a trivial solution is. Observe that in some of the examples discussed before it was very clear what trivial solutions look like. For example, trivial solutions to $k$-AP's are given by $x_1=\dots=x_k$, which any non-empty set would contain. In order to study the threshold we must avoid these kind of degenerate cases and understand what it means for the general setting.

For an admissible matrix $M\in \mM_{r,m}(\mathbb Z)$, consider the system $M\cdot \textbf{x}=\mathbf{0}$ and associate to each variable $x_i$ its corresponding index $i$. Let $\fp$ be a set partition of $\{1,\dots,m\}$ into $|\fp|$ blocks. Observe that $\fp$ defines a new system of equations $M_\fp\cdot \textbf{x'}=\mathbf{0}$, $M_{\fp}\in \mM_{r,|\fp|}(\mathbb Z)$, obtained after taking the original system $M\cdot \textbf{x}=\mathbf{0}$ and combining the variables of each block of $\fp$, i.e. summing up all columns related to the same block. We say that this new system of equations $M_\fp\cdot \textbf{x'}=\mathbf{0}$ is \emph{associated to $\fp$} and \emph{derived from $M\cdot \bx = \mathbf{0}$}.

A system associated to a partition $\fp$ encodes certain solutions of the original system with $m-|\fp|$ repeated coordinates. For a given solution $\textbf{x}$ of the system $M\cdot \textbf{x}=0$, we denote by $\fp(\textbf{x})$ the corresponding set partition of the indices $\{1,\dots,m\}$. In particular, if $\textbf{x}$ is a proper solution, then $\fp(\textbf{x})=\{\{1\},\dots,\{m\}\}$ and therefore $|\fp(\textbf{x})|=m$.

\SP

Observe that not every possible partition $\fp$ will come from a solution $\bx$ and not every system associated to a partition will have proper solutions. For example, if one considers the equation $x_1+x_2=x_3+x_4$ it is clear that the related partition $\fp=\{\{1\},\{2,3\},\{4\}\}$ (that is $x_2=x_3$) necessarily implies $x_1=x_4$, and thus the associated system will no longer be admissible (since it is neither irredundant nor non-degenerate). This observation is crucial in order to define what a trivial solution will be.

\begin{definition}\label{def:trivial}
We say that $\bx$ is a trivial solution of $M\cdot \bx=\mathbf{0}$ if
$$\rank{M_{\fp(\bx)}} < r.$$
We denote the set of all partitions stemming from some non-trivial solution by
$$\fP(M) = \left\{ \fp(\bx):\, M \cdot \bx = \mathbf{0}, \: \bx \text{ non-trivial} \right\}.$$
\end{definition}

Note that by definition $M_\fp$ is admissible for all $\fp \in \fP(M)$. Roughly speaking, our definition requires that the systems associated to our non-trivial solutions do not lose in complexity compared to the original system. Otherwise those solutions might in fact start appearing for smaller values of $p$. We already observed that trivial $k$-APs consisting of a single element would occur in any non-empty set.

Observe first that Definition~\ref{def:trivial} generalizes the notion of trivial solutions in the case $r=1$ introduced by Ruzsa in~\cite{Rz93}. Let us remark that previous generalizations have been made, like the one given by Shapira in~\cite{Sh06} for density regular systems of equations. Indeed, his definition of trivial solution is more restrictive: in his context $\bx$ is a trivial solution only if the system associated to $\mathfrak p(\bx)$ has zero rank (that is $\sum_{j\in P} a_{ij}=0$ for every $i$ and every $P\in \mathfrak p(\bx)$). As we will discuss later, this includes redundant solutions in some examples, but this does not affect his argument since he is interested in lower bounds for large sets avoiding solutions.

Let us discuss some examples to motivate our definition, i.e. to show that we are in fact no longer dealing with the same arithmetic structures when considering systems associated to trivial solutions. In Sidon sets, which are defined by the equation $x_1+x_2=x_3+x_4$, we can derive systems like $x_1+x_2=2x_3$ (namely $x_3=x_4$ in the original system) that give rise to non-trivial solutions since the rank of the associated system is still $1$. However, as said before, if one considers the partition $x_1=x_3,\ x_2=x_4$ then the associated system has rank $0$ and thus all solutions of this kind are trivial, which is consistent with the classical definition.

Next, let us discuss what trivial solutions look like for $B_h[g]$-sets. Recall that a set $\mA$ is no longer a $B_h[g]$-set if there exist $g+1$ (essentially different) representations of the same element as sums of $h$ elements of $\mA$. That is, there are elements $a_i\in \mA$, $1 \leq i \leq h(g+1)$ satisfying
$$a_1+\dots+a_h=a_{h+1}+\dots+a_{2h}=\dots=a_{hg+1}+\dots +a_{h(g+1)},$$
and all representations are pairwise different, so none of them are obtained after permuting the $h$ elements of another representation. Let us focus on $B_3[2]$ sets to illustrate what situations can occur. Here, we must avoid solutions to
$$x_1+x_2+x_3=x_4+x_5+x_6=x_7+x_8+x_9,$$
and we are excluding situations like for example $x_1=x_4,\ x_2=x_5,\ x_3=x_6$, with associated partition $\fp_1=\{\{1,4\},\{2,5\},\{3,6\},\{7\},\{8\},\{9\}\}$, since
$$M_{B_3[2]}=\left(\begin{smallmatrix}
      1&1&1&-1&-1&-1&0&0&0\\
0&0&0&1&1&1&-1&-1&-1
\end{smallmatrix}\right)\stackrel{x_1=x_4, x_2=x_5, x_3=x_6}{
\longrightarrow } \left(\begin{smallmatrix}
      0&0&0&0&0&0\\
1&1&1&-1&-1&-1
\end{smallmatrix}\right)=\left(M_{B_3[2]}\right)_{\fp_1}.$$
But we should not exclude for example solutions $x_1=x_3=x_5=x_7,\ x_4=x_8$, with partition $\fp_2=\{\{1,3,5,7\},\{2\},\{4,8\},\{6\},\{9\}\}$, which is still a valid solution since
$$M_{B_3[2]}=\left(\begin{smallmatrix}
      1&1&1&-1&-1&-1&0&0&0\\
0&0&0&1&1&1&-1&-1&-1
\end{smallmatrix}\right)\stackrel{x_1=x_3=x_5=x_7}{
\longrightarrow } \left(\begin{smallmatrix}
      1&1&-1&-1&0&0\\
0&0&1&1&-1&-1&
\end{smallmatrix}\right)\stackrel{x_4=x_8}{
\longrightarrow } \left(\begin{smallmatrix}
      1&1&-1&-1&0\\
0&0&0&1&-1
\end{smallmatrix}\right).$$
As we have seen in the Sidon case, different representations cannot have elements in common but the same representation can have repeated elements. If $h\geq 3$, we can also consider representations that have some elements in common but not all at once. As we observed before, the definition of trivial solutions in~\cite{Sh06} situates these two examples at the same level and clearly they should be considered different for our purposes.

\SP

Considering this definition of (non-)trivial solutions and the previous Lemma~\ref{fundamental}, one can already observe that the main contribution in our analysis of the threshold will come from proper solutions. The number of such solutions is, nevertheless, easier to count than the number of solutions with repeated coordinates (as we are dealing with general systems). The main difficulty will be to prove that the contribution of non-trivial solutions with repeated coordinates is negligible with respect to the total number of non-trivial solutions.

\section{Induced submatrices and related definitions}\label{submat}

We start this section by motivating the need for a definition of induced submatrices (not to be confused with the matrices associated to some partition, discussed on the previous section). Extending the results obtained for simple systems like that of Sidon sets~\cite{GJLR99}, one might expect the exponent of the threshold function of a given admissible system $M\cdot \textbf{x}=\mathbf{0}$ to be determined by the quotient of the number of variables of $M$ over the degrees of freedom in the system, which we will call the \emph{average degree} of $M$. However this exponent might not necessarily hold, as demonstrated by the example
\begin{equation}\label{eq:submatrix_example}
	M_1 = \left(\begin{array}{ccccccccc}
         1 & 1 & -1 & -1 & 0 & 0 & 0 & 0 & 0\\
         1 & 1 & 1 & 1 & 1 & 1 & 1 & 1 & -6
      \end{array}\right)\in \mM_{2,9}(\ZZ).
\end{equation}
We note that this system is of rank $2$ and has $9$ variables. Following the previous intuition, one might expect the exponent of the threshold function to be $-7/9$. The first row however implies that a solution to this system also fulfills the Sidon property ($x_1+x_2=x_3+x_4$) for which the stronger exponent $-3/4$ is known to hold. It follows that we should also consider the underlying information coming from so-called \textit{induced submatrices} in order to find an exponent that maximizes the average degree (see Definition~\ref{def:submatrix}).

\SP

In the general context, the induced submatrices are less clearly presented than in example~\eqref{eq:submatrix_example}. Intuitively, for any selection of rows one would try to set as many columns of these rows to zero through Gaussian elimination. Equivalently, one could fix some columns $Q\subseteq\{1,\dots, m\}$ and maximize the number of rows whose entries in the columns $\widebar{Q}=\{1,\dots, m\}\backslash Q$ can be set to zero through elimination such that one can disregard the corresponding variables.

For a proper definition, let $c_1, c_2, \dots, c_m$ be the columns of $M$. We have already introduced $M^{\widebar{Q}}$ as the matrix with columns $c_i$ for $i\in \{1,\dots,m\} \backslash Q$ so that $\rank{M^{\widebar{Q}}} = r - r_Q$. Denote by $b_1^Q, b_2^Q, \dots, b_r^Q$ the rows of $M^Q$ and by $b_1^{\widebar{Q}} , b_2^{\widebar{Q}} , \dots, b_r^{\widebar{Q}}$  the rows of $M^{\widebar{Q}}$. Assume without loss of generality that the first $r - r_Q$ rows of $M^{\widebar{Q}}$ are linearly independent. It follows that there exist $\delta_1^i, \delta_2^i, \dots, \delta_{r - r_Q}^i$ such that
\[ b_i^{\widebar{Q}} = \sum_{j=1}^{r - r_Q} \delta_j^i b_j^{\widebar{Q}} \  \text{ for }\  r - r_Q< i\leq r.\]
We use this notation for the following definition of our induced submatrix which is due to R\"odl and Ruci\'nski~\cite[Definition~7.1]{RR97}. Note that we slightly deviate from their notation.

\begin{definition}\label{def:submatrix}
For a given $\emptyset\neq Q\subseteq\{1,\dots, m\}$ s.t. $r_Q > 0$ we refer to the matrix defined by the rows
$$b_i^Q - \sum_{j=1}^{r - r_Q} \delta_j^i b_j^Q \  \text{ for }\  r - r_Q< i\leq r.$$
the induced submatrix $\Sub{M}{Q} \in \mM_{r_Q,|Q|}$ of $M$.
\end{definition}
Note that $\Sub{M}{Q}$ has $r_Q$ rows and $|Q|$ columns and that it is admissible. Also note that for $Q=\{1,\dots, m\}$ we have $\Sub{M}{Q} = M$, i.e. we obtain the full matrix. Let us illustrate with the following diagram how to obtain $\Sub{M}{Q}$ from $M$.
\begin{equation}\label{eq:diagram}
M \ \xrightarrow[\text{columns}]{\text{reordering}}\
\left(M^Q\ \Big|\ M^{\widebar{Q}}\,\right)
\xrightarrow[\text{of rows in }M^{\widebar{Q}}]{\text{gauss. elimination}}
\underset{=\delta(M)}{\underbrace{\left(\begin{array}{c|c} \cdots &\cdots\\ \hline
\Sub{M}{Q}& \textcolor{white}{\Big|}0
\end{array}\right)}}\begin{array}{l}\big]\ r-r_Q\text{ rows}\\ \textcolor{white}{\big|}\end{array}
\end{equation}

The following simple proposition now gives us the desired connection between a non-trivial solution $\bx = \left( x_1, \dots, x_m \right)$ of $M\cdot \bx=\textbf{0}$ and solutions to any induced submatrix $\Sub{M}{Q}$ if we set $\bx _Q = \left( x_i \right) _{i \in Q}$. Note that this is a slightly adapted version of \cite[Proposition~7.1]{RR97}.

\begin{proposition}\label{prop:connectionsubmatrix} Let $r<m$ and $M\in \mM_{r,m}(\mathbb Z)$ be an admissible matrix. For every $\emptyset\neq Q\subseteq\{1,\dots, m\}$ and non-trivial solution $\bx$ of the system $M\cdot \bx=\textbf{0}$ we have that $\bx _Q$ is a non-trivial solution of $\Sub{M}{Q} \cdot \bx _Q = \textbf{0}$.
\end{proposition}

\begin{proof}
Let $\bx$ be a non-trivial solution of of the system $M\cdot \bx=\textbf{0}$. Without loss of generality we assume that the columns indexed by $Q$ are the first $|Q|$ columns of $M$ and that the indices of the $r-r_Q$ linearly independent rows of $M^{\widebar{Q}}$ are also the indices of the first $r-r_Q$ rows. Now let $\delta (M)$ be the matrix obtained by applying the linear transformations described in Definition~\ref{def:submatrix} to the whole matrix $M$ (see~\eqref{eq:diagram} for a visual definition of $\delta(M)$).

We note that by construction the matrix $\Sub{M}{Q}$ can be found in the lower left of $\delta (M)$. Now since we only applied linear transformations to the rows of $M$ in order to obtain $\delta (M)$ we still have $\delta(M)\cdot \bx=\textbf{0}$ and therefore also $\Sub{M}{Q} \cdot \bx _Q = \textbf{0}$.

Next we observe that $\delta \left( M_{\fp(\bx)} \right) =  \delta (M) _{\fp(\bx)}$ and that $M_{\fp(\bx)}$ has rank $r$ since $\bx$ was a non-trivial solution. Noting again that we are using only linear transformations we can state that $\delta \left( M_{\fp(\bx)} \right)$ also has rank $r$ and hence $\delta (M) _{\fp(\bx)}$ does as well. Now since $\delta (M) _{\fp(\bx)}$ is of maximum possible rank it follows that $\Sub{M}{Q}_{\fp(\bx)}$ has the same rank as $\Sub{M}{Q}$. We can conclude that $\bx _Q$ is indeed a non-trivial solution of $\Sub{M}{Q} \cdot \bx _Q = \textbf{0}$.
\end{proof}
The main implication of this proposition will be that if for some induced submatrix $\Sub{M}{Q}$ there do not exist non-trivial solutions then there will not be non-trivial solutions for $M\cdot \bx=\textbf{0}$ either. We observe that in Definition~\ref{def:maxparameter} of $c(M)$ we are maximizing the average degree (ratio between number of variables and degrees of freedom) over all induced submatrices. Considering the previous proposition this should be no surprise and in fact is the core idea of the proof of Theorem~\ref{thm:threshold}.

\SP

Before continuing, we will introduce the definition of a (strictly) balanced system that is needed in order to state Theorem~\ref{thm:local}. The definition basically requires that the average degree  of any induced submatrix of a given derived system is (strictly) smaller than that of the original system. In particular, this implies that $c(M) = m/(m-r)$.

\begin{definition}\label{def:balanced}
We call an admissible matrix $M \in \mM_{r,m}(\ZZ)$ balanced if for every $\fp \in \fP(M)$ with $2 \leq |Q| < | \fp |$ we have the inequality
$$\frac{|Q|}{|Q|-\rank{\Sub{M_\fp}{Q}}} \leq \frac{m}{m-r}.$$
If the inequality is strict, we say that the system is strictly balanced.
\end{definition}

This is reminiscent of a similar definition obtained in graph theory~\cite{ER60}. Note that it differs slightly from R\"odl - Ruci\'nski definition of strictly balanced systems~\cite[Definition~7.2]{RR97}. There exist simple examples showing that this is a real difference and in fact it will be crucial in the proof of Theorem~\ref{thm:local}.

As an example, consider the homogeneous system given by the matrix
$$M_2 = \left(\begin{array}{cccccccc}
         1 & -2 & 1 & 0 & 0 & 0 & 0 & 0\\
         0 & 1 & 1 & -1 & -1 & 0 & 0 & 0\\
         1 & 1 & 1 & 1 & 1 & 1 & 1 & -7
      \end{array}\right)\in \mM_{3,8}(\ZZ).$$
We note that the full system would give us the exponent $-5/8$. Observe that the system implicitly contains the 3-AP condition as an induced submatrix (if we choose $Q_\text{3-AP}=\left\{ 1,2,3 \right\}$) as well as the Sidon condition (if we select $Q_\text{Sidon}=\left\{ 2,3,4,5 \right\}$). The exponents obtained for either induced submatrix are lower than that of the full system. However if we choose $Q=\left\{ 1,2,3,4,5 \right\}$ we now have $r_Q = r-1$ and hence
$$\Sub{M_2}{Q} = \left(\begin{array}{ccccc}
		1 & -2 & 1 & 0 & 0\\
		0 & 1 & 1 & -1 & -1
	\end{array}\right)\in \mM_{2,5}(\ZZ).$$
This gives us the larger exponent $-3/5$, which means that the full system is not balanced. It is in fact the largest exponent we can achieve, so in this case we have $c(M_2) = 5/3$.

\section{Proof of Theorem~\ref{thm:threshold} }\label{sec:thresholdproof}

\subsection{The 0-statement.} \label{0statement}

In order to prove the 0-statement in Theorem~\ref{thm:threshold} we note that by Definition~\ref{def:maxparameter} one can pick a specific set of columns $Q \subseteq \{1,\dots, m\}$ satisfying $c(M) = |Q|/(|Q| - r_Q)$. For simplicity we write $q = |Q|$ so that $r_Q = (1 - 1/c(M))q$ and we also denote the corresponding induced submatrix $\Sub{M}{Q} \in \mM_{r_Q,q}(\mathbb Z)$ by $B$.

The proof proceeds in the following way: we will apply the First Moment Method in order to show that the probability of $\mA^q$ containing a non-trivial solution to $B\cdot \bx=\textbf{0}$ tends to zero for random sets $\mA\subseteq\left[n\right]$ where each element is chosen independently with probability $p = o\left(n^{-1/c(M)}\right)$. By Proposition~\ref{prop:connectionsubmatrix} it will follow that the probability of $\mA^m$ containing a solution to our original system $M\cdot \bx=\textbf{0}$ also goes to zero. That is, the random variable which counts the number of non-trivial solutions goes to $0$ asymptotically almost surely for $p = o\left(n^{-1/c(M)}\right)$.

Let $S_B$ be the set of (essentially) different non-trivial solutions in $\mA^q$ (possibly with repeated coordinates) of $B\cdot \textbf{x}=0$ and denote the corresponding counting variable by $\bX_B$. Note that $S_B$ is the set of solutions considered as vectors modulo certain symmetries arising from $M$ (see Section~\ref{sec:countingsolutions} for a detailed description of the set of symmetries defining this equivalence relation). We can write the counting random variable as
\begin{equation}\label{EE}
\bX_B=\sum_{\bx\in S_B}\II_{\bx},
\end{equation}
where $\II_{\bx}$ denotes the indicator random variable for the event $\bx \in \mA^{q}$ which we denote as $E_{\bx}$.

It is clear that $\EE{\II_\bx}=\PP{E_{\bx}}=p^{|\fp(\textbf{\bx})|}$, where $\fp(\textbf{x})$ is the set partition associated to $\bx$ (see Section~\ref{trivialsol}). Denote by $S_{B_\fp}^\star$ the set of (essentially) different proper solutions of some system $B_\fp \cdot \bx = \bf{0}$ where $\fp \in \fP(B)$. Recall that $|\fp|\ge \rank{B}$ for every $\fp\in\fP(B)$, since it follows from the definition of trivial solution that the derived system $B_{\fp}$ must not decrease in rank. Now by splitting up the sum in \eqref{EE} by the size of the corresponding set partition we have
\begin{equation}\label{eq:sum-ind}
\EE{\bX_B}
= \sum_{\bx \in S_B}\PP{E_{\bx}}
= \sum_{s=\rank{B}}^q \enspace \sum_{\substack{\fp \in \fP(B) \\ |\fp| = s}} \enspace \sum_{\substack{ \bx\in S_B \\ \fp(\textbf{x}) = \fp }}p^s
= \sum_{s=\rank{B}}^q \enspace \sum_{\substack{\fp \in \fP(B) \\ |\fp| = s}} \enspace |S_{B_\fp}^\star| p^s.
\end{equation}
Since we are considering non-trivial solutions and we have assumed that $q/(q-\rank{B}) = c(M)$, we have $\rank{B_\fp} = \rank{B} > (1 - 1/c(M)) q \geq (1 - 1/c(M)) |\fp|$ for $\fp \in \fP(B)$ and can apply Equation~\eqref{eq:trivialupperbound} to obtain the estimate
$$\EE{\bX_B} = \sum_{s=q-q/c(M)}^q O \left( n^{ -s/c(M)} p^s \right) = \sum_{s=q-q/c(M)}^q O \left( \frac{p}{n^{1/c(M)}} \right)^s \to 0.$$

\subsection{The 1-statement.} \label{1statement}

Instead of dealing with the random variable $\bX$ which counts the number of non-trivial (essentially different) solutions to $M\cdot \bx = \mathbf{0}$, we simplify the problem by only focusing on proper solutions. Denote the set of (essentialy different) proper solutions by $S^\star_M$ and write $\bX^\star$ for the random variable counting the number of such solutions. We will use the Second Moment Method described in Subsection~\ref{sub:2nd-moment} to prove that $\bX^\star>0$ asymptotically almost surely for $p=\omega(n^{-1/c(M)})$. This clearly also implies that $\bX>0$.

We start by showing that the first moment of $\bX^\star$ goes to infinity. Note that by simply choosing $Q=\{1,\dots, m\}$ in Definition~\ref{def:maxparameter} we obtain $c(M) \geq m/(m-r)$. Therefore $p=\omega\big(n^{-1/c(M)}\big)$ also implies $p=\omega(n^{-(m-r)/m})$ and we can apply Lemma~\ref{fundamental} to get
\begin{align}\label{eq:X_0}
	\EE{\bX^\star} & = \sum_{\bx \in S^\star_M} \PP{E_{\bx}} = \frac{\Vol{\mP_M}}{\sigma(M)}n^{m-r}p^{m}(1+o(1)) \nonumber \\
	& = \frac{\Vol{\mP_M}}{\sigma(M)} \left( \frac{p}{n^{-(m-r)/m}} \right)^m (1+o(1)) \to \infty.
\end{align}

Now we must carefully study the second moment of the variable $\bX^\star$ in order to conclude that in fact $\bX^\star>0$ asymptotically almost surely for $p=\omega\big(n^{-1/c(M)}\big)$. We know that the events $E_{\bx}$ and $E_{\by}$ of two different proper solutions $\bx,\,\by \in S^\star_M$ are not independent if and only if $0<|\{x_1,\dots,x_m\}\cap \{y_1,\dots, y_m \}|<m$. That is some (but not all) of their coordinates coincide. Following the notation of Subsection~\ref{sub:2nd-moment} we write $\bx \sim \by$ for two such solutions and note that the number of coincidences must be at least $1$ to guarantee that $\bx\sim \by$ and at most $m-1$ so that $\bx\ne\by$. As defined in Equation~\eqref{eq:delta}, we will have to study the quantity
$$\Delta^\star = \sum_{\substack{\bx,\by \in S^\star_M \\ \bx\sim \by}}\PP{E_{\bx}\wedge E_{\by}},$$
where the sum is taken over all intersecting pairs $\bx,\by\in S^\star_M$.

In order to show that $\Delta^\star = o \big( \EE{\bX^\star}^2 \big)$, fix a solution $\bx \in S^\star_M$ as well as a non-empty bipartite matching $G=(\{1,\dots, m\},\{1,\dots, m\},E)$ where $1 < |E| < m$. Now consider all solutions $\by \in S^\star_M$ whose coincidences with $\bx$ are indicated by $G$, i.e. $\bx_i = \by_j$ if and only if $\overline{ij} \in E$. It is clear that finding such solutions is equivalent to finding proper solutions to the linear system
$$M^{\widebar{Q(G)}} \cdot \by_{\widebar{Q(G)}} = \mathbf{b},$$
where $Q(G) = \left\{ j \mid (i,j) \in E \text{ for some } j \in \{1,\dots, m\} \right\}$ and $\mathbf{b} = -\sum_{(i,j) \in E} \bx_i \cdot c_j$ where $c_j$ denotes the columns of $M$ as before. We note that, as previously defined, the rank of $M^{\widebar{Q(G)}}$ is $r_{\widebar{Q(G)}}$. Next we observe that by Equation~\eqref{eq:trivialupperbound} the number of such solutions is bounded by
$$O\big(n^{|\widebar{Q(G)}|-r_{\widebar{Q(G)}}}\big) = \big(n^{m-|Q(G)|-r_{\widebar{Q(G)}}}\big)$$
and that $\PP{E_{\bx}\wedge E_{\by}} = p^{2m-|E|} = p^{2m-|Q(G)|}$. For both of these values only the set of columns $Q(G)$ as well as its cardinality are of importance and not the exact bipartite graph $G$ indicating the coincidences. We therefore simply write $\bx \sim_Q \by$ whenever the coincidences between two solutions $\bx, \by \in S^\star_M$ are given by some non-empty bipartite graph $G$ satisfying $Q(G) = Q$.  Noting that by Lemma~\ref{fundamental} the number of proper solutions $\bx \in S^\star_M$ is $\Theta(n^{m-r})$, we can now write
\begin{align}\label{eq:final_brun}
\Delta^\star & = \sum_{\substack{\bx,\by \in S^\star_M \\ \bx\sim \by}}\PP{E_{\bx}\wedge E_{\by}} = \sum_{\bx \in S^\star_M} \sum_{\emptyset \neq Q \subsetneq [m]} \sum_{\substack{\by \in S^\star_M \\ \bx \sim_Q \by}} p^{2m-|Q|} \nonumber \\
& = \sum_{\emptyset \neq Q \subsetneq [m]} \sum_{\bx \in S^\star_M} O\left(n^{m-|Q|-(r-r_Q)}\right) p^{2m-|Q|} = \sum_{\emptyset \neq Q \subsetneq [m]} O\left(n^{2m-r-|Q|-(r-r_Q)}\right) p^{2m-|Q|} \nonumber \\
& = O\left( \left( n^{m-r}p^m \right)^2 \right) \sum_{\emptyset \neq Q \subsetneq [m]} \left( \frac{n^{-(|Q|-r_Q)/|Q|}}{p} \right)^{|Q|} = o \left( \EE{\bX^\star}^2 \right),
\end{align}
since $p = \omega(n^{-1/c(M)})$ also implies $p = \omega(n^{-(|Q|-r_Q)/|Q|})$ for all $Q \subseteq \{1,\dots, m\}$ by Definition~\ref{def:maxparameter} and $\EE{\bX^\star} = \Theta (n^{m-r}p^m)$ by Equation~\eqref{eq:X_0}. By the Second Moment Method it follows from~\eqref{eq:final_brun}  that $\bX^\star\sim\EE{\bX^\star}$ asymptotically almost surely. In particular we have that, in this range, $\bX>\bX^\star>0$ with probability tending to $1$.

\qed

\section{Proof of Theorem~\ref{thm:local}}\label{Brun-proof}

\subsection{Sufficiency}

We will apply Brun's Sieve (see Theorem~\ref{Thm:Brun}) in order to show the statement. Note that we require our system to be strictly balanced so in Definition~\ref{def:balanced} we can pick $\fp = \{1,\dots, m\}$ (the partition of a proper solution) and in particular we have that $|Q|/(|Q|-r_Q) < m/(m-r)$ for all $Q \subsetneq \{1,\dots, m\}$ and hence $c(M) = m/(m-r)$.

We start by developing the first moment of the random variable $\bX$ which counts the number of essentially different non-trivial solutions to $M\cdot \bx = \mathbf{0}$ in $\mA^m$ for $p = Cn^{-1/c(M)} = Cn^{-(m-r)/m}$. As before we split it up by the number of repeated coordinates in the solutions and apply Lemma~\ref{fundamental}:
\begin{align*}
\EE{\bX} &= \sum_{\bx \in S_M} \PP{E_{\bx}} = \sum_{s=1}^m \enspace \sum_{\substack{\fp \in \fP(M) \\ |\fp| = s}} \enspace \sum_{\bx\in S_{M_\fp}^\star} p^s\\
&= \tfrac{\Vol{\mP_M}}{\sigma(M)} n^{m-r}p^{m}(1+o(1)) + \sum_{s=1}^{m-1} \enspace \sum_{\substack{\fp \in \fP(M) \\ |\fp| = s}} o\left( n^{s - \tfrac{r}{m}s} \right) p^s\nonumber \\
& = C^m \tfrac{\Vol{\mP_M}}{\sigma(M)}(1+o(1)) +o(1) = \mu(1+o(1))
\end{align*}
where we set $\mu = C^m \, \Vol{\mP_M}/\sigma(M)$. Note that in the third equality the first part refers to proper solutions and the second part to non-trivial solutions $\bx$ with repeated coordinates for which we have $\rank{M_{\fp (\bx)}} = r > r \, |\fp (\bx)| / m$ by Definition~\ref{def:trivial}.

Let $\EE{\bX}_t$ denote the expected number of ordered $t$-tuples of solutions and split it up into three parts:
$$\EE{\bX}_t = \EE{\bX}_t' + \EE{\bX}_t'' + \EE{\bX}_t'''.$$
The first part $\EE{\bX}_t'$ refers to $t$-tuples of pairwise disjoint proper solutions, $\EE{\bX}_t''$ refers to $t$-tuples of pairwise disjoint non-trivial solutions of which at least one is not proper and $\EE{\bX}_t'''$ refers to $t$-tuples of non-trivial solutions in which at least two share a coordinate. We will compute each of them to show that $\EE{\bX}_t = \mu^t \left(1 + o(1)\right)$ from which we can follow that $S^{(t)} \to \mu^t/t!$ so that Theorem~\ref{Thm:Brun} applies.

In order to compute the expected number of ordered $t$-tuples of solutions we introduce the notion of \emph{compounded} systems of linear equations. For this consider two matrices $A$ and $B$ with $m_A$ and $m_B$ columns respectively as well as an incomplete bipartite matching $G=([m_A],[m_B],E)$ where $E = \left\{(i_1,j_1),\dots,(i_k,j_k)\right\}$ for some $k<\min(m_A,m_B)$ (note the similarity to the previous proof). Denote the columns of $A$ and $B$ by $a_i$ and $b_j$. Write also $A^{[m_A] \backslash \left\{ i_{1},\ldots,i_{k}\right\}}$ and $B^{[m_B] \backslash \left\{ j_{1},\ldots,j_{k}\right\}}$ for the matrices obtained from $A$ and $B$ by removing those columns indexed by $\left\{ i_{1},\ldots,i_{k}\right\}$ or $\left\{ j_{1},\ldots,j_{k}\right\}$ respectively.  Define now the following matrix:
$$A \times_G B = \left(\begin{array}{ccccc}
A^{[m_A] \backslash \left\{ i_{1},\ldots,i_{k}\right\} } & a_{i_{1}} & \ldots & a_{i_{k}} & 0\\
0 & b_{j_{1}} & \ldots & b_{j_{k}} & B^{[m_B] \backslash \left\{ j_{1},\ldots,j_{k}\right\} }
\end{array}\right).$$
For our application the actual matching and hence the concrete type of overlap will be irrelevant. What matters is whether the systems are disjoint or not, so if the bipartite graph is empty or if there is some actual overlap. We therefore simply omit the graph in our notation and write $A\times B$ when compounding matrices and $A\, \dot{\times}\, B$ to specify when they are being compounded without overlap. Note that this operator is not commutative or associative (and in fact strongly depends on the size of the respective matrices) and we will write $A \times B \times C = \left( A \times B \right) \times C$

Using this we note that finding $t$-tuples of pairwise disjoint proper solutions to $M$ is equivalent to finding proper solutions of the compounded system $M \dot{\times} \overset{t}{\ldots} \dot{\times} M$. The resulting system in $tm$ variables $M \dot{\times} \overset{t}{\ldots} \dot{\times} M \cdot \bx = \mathbf{0}$ is trivially admissible and has rank $tr$. It is also easy to see that $\sigma(M \dot{\times} \overset{t}{\ldots} \dot{\times} M) = \sigma(M)^t$  and $\Vol{\mP_{M \dot{\times} \overset{t}{\ldots} \dot{\times} M}} = \Vol{P_M}^t$ and therefore we can apply Lemma~\ref{fundamental} in order to get
$$\EE{\bX}_t' = \tfrac{\Vol{\mP_{M \dot{\times} \overset{t}{\ldots} \dot{\times} M}}}{\sigma(M \dot{\times} \overset{t}{\ldots} \dot{\times} M)} n^{tm-tr} p^{tm} \left(1 + o(1) \right) = \mu^t \left(1 + o(1) \right).$$

Next consider all $t$-tuples of pairwise disjoint non-trivial solutions of which at least one is not proper. This means we are considering all compounded systems $M_{\fp_1} \dot{\times} \ldots \dot{\times} M_{\fp_t}$ where $\fp_i \in \fP(M)$ and at least one of them is not equal to $\{1,\dots, m\}$, i.e. it does not come from a proper solution. We note that the compounded system trivially is non-degenerate. The system also is irredundant since otherwise it would not actually be associated to some $t$-tuples of pairwise disjoint non-trivial solutions and hence would not be part of $\EE{\bX}_t''$, so in particular it is admissible. Its number of columns is $\sum_{i} |\fp_i|$ and – since the partitions come from non-trivial solutions – its rank is $\sum_{i} \rank{M_{\fp_i}} = tr$. It follows by Lemma~\ref{fundamental} that
\begin{align*}
\EE{\bX}_t'' & = O \big( n^{\sum_{i} |\fp_i| - tr}p^{\sum_{i} |\fp_i|} \big)  = O \big( n^{tr(\sum_{i} |\fp_i|/tm - 1 )} \big)  \to 0
\end{align*}
where the limit follows from the fact that $\sum_{i} |\fp_i| < tm$. Here we have used the assumption that one of the partitions does not come from a proper solution, so there exists a partition $|\fp_i|$ which satisfies $|\fp_i| < m$.

So far we have not used the condition that our system is strictly balanced. This requirement will be of importance now when considering $t$-tuples of non-trivial solution in which at least two solutions share a coordinate. This means we are considering all compounded systems $M_{\fp_1} \times \ldots \times M_{\fp_t}$ where $\fp_i$ are partitions of the columns resulting from some non-trivial solutions and for at least one of the compound operators the implied bipartite matching is not empty. This means at least one of the systems $M_{\fp_i}$ overlaps with the rest of the construction in some way. The irredundancy of the system follows trivially as in the previous case.

Let us consider in general terms what happens when we compound two admissible systems $A \in \mM_{r_A,m_A}(\mathbb Z)$ and $B \in \mM_{r_B,m_B}(\mathbb Z)$. We have already established that $A \dot{\times} B$ is a system in $m_A+m_B$ variables with rank $r_A+r_B$. Now let us assume there is some overlap, that is some columns $\emptyset \neq Q \subsetneq [m_B]$ of $B$ are matched to columns of $A$ in the implied bipartite matching. The compounded system will have $m_a + m_b - |Q|$ columns and be of rank at least $r_A + r_B - r_Q(B)$. This follows easily since if some rows in the compounded system stemming from $B$ are linearly dependent, then first their coordinates in $\widebar{Q}$ have to be linearly dependent on other rows stemming from $B$. Using the above notation gives the upper bound $r_Q(B)$ for the number of rows that can become linearly dependent by compounding the two systems.

We know that $M_{\fp_1} \times \ldots \times M_{\fp_t}$ is a system of equations in $tm-\beta$ variables of rank $tr-\alpha$ for some $\alpha,\beta \in \NN_0$. Since we are only considering systems with some overlap we assume $\beta > 0$. Using the previous observations we can induce over $1 \leq i \leq t$ to show that $\beta \, r/m > \alpha$. Assume w.l.o.g. that the first two systems $M_{\fp_1}$ and $M_{\fp_2}$ overlap in some columns $\emptyset \neq  Q \subsetneq [|\fp_2|]$ of $M_{\fp_2}$. It follows that $M_{\fp_1} \times M_{\fp_2}$ has $|\fp_1| + |\fp_2| - |Q|$ columns and by the previous observation is of rank at least $2r - r_Q(M_{\fp_2})$. Since our system is strictly balanced we have by Definition~\ref{def:balanced} that
$$\frac{|Q|}{|Q| - r_Q(M_{\fp_2})} < \frac{m}{m - r} \quad  \Leftrightarrow \quad r_Q(M_{\fp_2}) < \tfrac{r}{m}|Q|.$$
This means that the rank of $M_{\fp_1} \times M_{\fp_2}$ is strictly greater than $2r - |Q| \, r/m$. Therefore the first step of the induction is complete. Next we assume that $M_{\fp_1} \times \ldots \times M_{\fp_k}$ has $km-\beta$ variables and is of rank $kr-\alpha$ for $1<k<t$ and some $\beta \, r/m > \alpha$, $\beta > 0$. We compound $M_{\fp_{k+1}}$ with $M_{\fp_1} \times \ldots \times M_{\fp_k}$ where the overlap is indicated by $Q \subsetneq [|\fp_{k+1}|]$ and $Q = \emptyset$ is possible. It follows by the same arguments as before that $M_{\fp_1} \times \ldots \times M_{\fp_{k+1}}$ has $(k+1)m - (\beta + |Q|)$ columns and is of rank strictly greater than
$$kr - \alpha + r-r_Q(M_{\fp_{k+1}}) \geq (k+1)r - (\alpha + \tfrac{r}{m}|Q|)$$
since $\fp_{k+1}$ again comes from a non-trivial solution. Obviously we still have
$$\tfrac{r}{m} (\beta + |Q|) > \alpha + \tfrac{r}{m}|Q|$$
and therefore the induction is complete. Now using the fact that $\beta \, r/m > \alpha$ we can apply Lemma~\ref{fundamental} and state that
\begin{align*}
\EE{\bX}_t''' & = O \left( n^{(tm - \beta) - (tr - \alpha)}p^{tm - \beta} \right) \\
& = O \left( n^{t(m-r) + (\alpha - \beta)}n^{t(r-m) - \beta(r/m-1)} \right) \\
& = O \left( n^{\alpha - \tfrac{r}{m}\beta}\right) \to 0.
\end{align*}
Taken together it follows that
$$\EE{\bX}_t = \EE{\bX}_t' + \EE{\bX}_t'' + \EE{\bX}_t''' = \mu^t \left(1 + o(1) \right)$$
and since $S^{(t)} = \EE{\bX}_t / t!$ we can apply Brun's sieve to deduce the statement of the theorem.

\subsection{Necessity}

It remains to show that the strictly balanced condition is in fact necessary. If the system is balanced but not strictly balanced, i.e. $c(M) = m/(m-r)$ but there exists an induced submatrix also attaining this value, we again split $\EE{\bX}_t$ into three parts
$$\EE{\bX}_t = \EE{\bX}_t' + \EE{\bX}_t'' + \EE{\bX}_t'''$$
as before. Observe that $\EE{\bX}_t' = \mu (1+o(1))$ and $\EE{\bX}_t'' = o(1)$ are unchanged from the previous computations since we did not rely on the strictly balanced condition for their computation. Further continuing the notation from the proof of sufficiency, we know that the compounded systems $M_{\fp_1} \times \ldots \times M_{\fp_t}$ considered in $\EE{\bX}_t'''$ have $tm-\beta$ variables and are of rank $tr-\alpha$ for some $\alpha,\beta \in \NN_0$ s.t. $\beta > 0$. Doing a simple induction as before we can show that $\beta \, r/m \geq \alpha$ since the system is balanced. Note that previously we had a strict inequality since we were dealing with strictly balanced systems. If the inequality is strict, the compounded system is negligible as before. However, since by assumption our system is balanced but not strictly balanced, there can also exist compounded systems for which we have equality, i.e. $\beta \, r/m = \alpha$. Note that there is obviously a bounded number of $\alpha$ and $\beta$ for each $t$ so there is a finite number of these compounded systems. These (if they exists) would each contribute in the order of a constant, since
$$\Theta \left( n^{(tm - \beta) - (tr - \alpha)}p^{tm - \beta} \right) = \Theta \left( n^{tm - tr + \alpha - \beta}n^{tr - \tfrac{r}{m}\beta - tm + \beta} \right) \\ = \Theta \left( n^{\alpha - \tfrac{r}{m}\beta}\right) = \Theta(1).$$
We would like to conclude from this, that there exists a positive constant $c_t > 0$ such that $\EE{\bX}_t''' = o(1) + c_t$. However, we first need to verify that these compounded systems actually occur, i.e. that they come from some tuple of solutions with the correct overlap. Assume a proper solution $\bx$ of $M \cdot \bx = \mathbf{0}$ and fix some $\emptyset \neq Q \subsetneq \{1,\dots, m\}$ coordinates belonging to an induced submatrix for which $|Q|/(|Q|-\rank{\Sub{M}{Q}}) = m/(m-r)$. We have $r_Q = |Q| \, r/m$ and therefore one can easily show that $(m-|Q|) - (r-r_Q) > 0$. This means there is at least a degree of freedom in the remaining un-fixed coordinates, so it is possible to have two solutions overlap exactly in $Q$ but not any other coordinates and the compounded system does actually occur in $\EE{\bX}_t'''$.

Combining the previous observations, it follows that
$$\EE{\bX}_t = \EE{\bX}_t' + \EE{\bX}_t'' + \EE{\bX}_t''' = \mu^t \left(1 + o(1) \right) + c_t.$$
Obviously, for each $n>0$, fixed $s$ and $0\leq t\leq s$ the values $\EE{\bX^t}$ are moments of a random variable. Hence, they satisfy Stieltjes condition (see \cite{RS65}), which is preserved by taking limits. Consequently, for each $t$ the sequence $\mu^t \left(1 + o(1) \right) + c_t$ converges to the $t$-th moment of a certain random variable. Due to Carleman's condition, this random variable is indeed uniquely determined. Finally, the limit of the sequence $\bX$ does not have a moment sequence equal to $\{\mu^t/t!\}_{t\geq 1}$, but it is determined by its moments. We conclude that we cannot have convergence in distribution towards a Poisson distributed random variable.

To conclude the analysis, the unbalanced case can be deduced by using a similar argument to the one on the balanced case by conveniently rescaling the random variable and showing that it does not converge in distribution to a Poisson random variable. The details are the same as in the proof of \cite[Theorem 5]{Ru90}.

\section{The computation of $\Vol{\mP_M}$}\label{c(M)}

In this section we consider the question of computing the constants $\Vol{\mP_M}$ involved in Theorem~\ref{thm:threshold}. As we have shown in previous sections, the constant $\Vol{\mP_M}$ is the volume of the polytope defined by the equations $M\cdot \bx =0$, where the coordinates of the vector $\bx$ belong to the closed interval $[0,1]$.

We study the $k$-sum free sets as a warm up. Note that the $k$-barycentric case could be treated with the same ideas. Secondly we analyze Ehrhart's Polynomial for the polytope associated to $k$-AP's by means of elementary arguments. For $B_{h}[g]$-sets, we obtain an exact formula by means of Vandermonde's determinants. Finally, the volume in the case of $k$-cubes is not analyzed here, but observe that the volume can be deduced in this case from the results of~\cite{Sa07}.

\subsection{$k$-sum-free sets.}

As a toy example, let us compute the volume of the polytope associated to sum-free sets, obtained from the linear equation $x_1+x_2=x_3$, $0\leq x_i \leq 1$. The associated polytope can be defined as follows
$$\mP_{1-SF}=\{(x_1,x_3)\,:\, 0\leq x_1\leq x_3\leq 1\}\subset \mathbb R^2,$$
since $x_2=x_3-x_1\in [0,1]$ for any $(x_1,x_3)\in \mP_{1-SF}$. Clearly $\mP_{1-SF}$ is an integral polytope, since it is in fact the triangle with vertices $(0,0),\ (0,1)$ and $(1,1)$, and an easy computation gives a volume equal to $1/2$.

However, let us obtain this value by means of interpolation arguments. It follows from Ehrhart's Theorem that $\left|n\cdot \mP_{1-SF} \cap\,\mathbb{Z}^{2}\right|=f(n)$ for a polynomial $f$ of degree $\text{dim}(\mP_{1-SF})=2$; namely $f(n)=\Vol{\mP_{1-SF}}n^2+bn+c$. It is clear that $f(0)=|\{(0,0)\}|=1$ (which gives $c=1$), $f(1)=f(0)+|\{(0,1),(1,1)\}|=3$ (thus $b=2-\Vol{\mP_{1-SF}}$) and $f(2)=f(1)+|\{(0,2),(1,2), (2,2)\}|=6$. Therefore
$$f(2)=4\Vol{\mP_{1-SF}}+2b+c=2\Vol{\mP_{1-SF}}+5=6\ \Longrightarrow\ \Vol{\mP_{1-SF}}=\frac{1}{2},$$
as we wanted to show.

The case $k>1$ is slightly different: here we consider the set
$$\mP_{k-SF}=\{(x_1,x_3)\,:\, 0\leq kx_3- x_1\leq 1,\, 0\leq x_1,x_3\leq 1\}\subset \mathbb R^2,$$
which is a parallelogram instead of a triangle. Its area is equal to $1/k$. The main difference is that in the first case we obtain a polynomial, despite in the second case we may obtain a \emph{pseudopolynomial}.

\SP

We continue computing $\Vol{\mP_M}$ in the case of $k$-AP's and also for $B_{h}[g]$ sets. In the first case by elementary means we obtain the closed expression for the volume. In the former case, we apply interpolation arguments to obtain a general expression in terms of determinants.

\subsection{$k$-AP free sets.}

This family has been studied widely. For instance, the following result is also implicitly stated in~\cite{Sa07}. For completeness we include the analysis here. As we have seen earlier, a $k$-AP
$$x_1=a,\,x_2=a+d,\,x_3=a+2d,\dots, x_{k}=a+(k-1)d$$
can be expressed as a solution $\bx=(x_1,\dots,x_k)$ of a linear system of rank $k-2$ in $k$ variables. We can count the number of $k$-AP with elements in $[n] \cup\{0\}$ by direct counting:

\begin{proposition}\label{k-AP}
For any integer $k\geq 3$ the number of $k$-AP (including trivial ones) in $[n]\cup\{0\}$ is given by
$$(n+1)\left(\left\lfloor\frac{n}{k-1}\right\rfloor+1\right)-\frac{k-1}{2}\left(\left\lfloor\frac{n}{k-1}\right\rfloor^2+\left\lfloor\frac{n}{k-1}\right\rfloor\right).$$
\end{proposition}

\begin{proof}
Observe that any $k$-AP is of the form $\{a,a+d,\dots,a+(k-1)d\}$ where $a\in [n]\cup\{0\}$ and $d\in \{0,1,2,\dots,\lfloor n/(k-1) \rfloor\}$, since
$$ 0,\, \lfloor\tfrac{n}{k-1}\rfloor,\,2\lfloor\tfrac{n}{k-1}\rfloor,\dots,\,(k-1)\lfloor\tfrac{n}{k-1}\rfloor\leq n$$
is a $k$-AP and the equality holds for multiples of $k-1$. Additionally, for a given $d$ we have that $$\{0, d,\dots, (k-1)d)\},\ \{1, 1+d,\dots, 1+(k-1)d)\},\ \dots,\ \{n-(k-1)d,n-(k-2)d\dots,n\}$$ are the only $k$-AP with common difference $d$. Thus the total number of $k$-AP is given by
$$\sum_{d=0}^{\lfloor\tfrac{n}{k-1}\rfloor} \left(n+1-(k-1)d\right)=(n+1)\left(\left\lfloor\frac{n}{k-1}\right\rfloor+1\right)-\frac{k-1}{2}\left(\left\lfloor\frac{n}{k-1} \right\rfloor^2+\left\lfloor\frac{n}{k-1}\right\rfloor\right).$$
\end{proof}

\begin{corollary}
The polytope associated to the $k$-AP condition has volume $1/(k-1)$.
\end{corollary}

\begin{proof}
Let $\mP_k$ denote the associated polytope. By Ehrhart's Theorem it follows that the number of $k$-AP in $[n]\cup\{0\}$ is equal to $\Vol{\mP_k} n^{2}+O(n)$. Note that using the number of solutions as counted in Proposition~\ref{k-AP} satisfies
$$(n+1)\left(\left\lfloor\frac{n}{k-1}\right\rfloor+1\right)-\frac{k-1}{2}\left(\left\lfloor\frac{n}{k-1}\right\rfloor^2+\left\lfloor\frac{n}{k-1}\right\rfloor\right) = \frac{1}{2(k-1)} n^2 + O(n).$$
However Proposition~\ref{k-AP} only counted solutions up to symmetry so we needed to take multiply by a factor of $\sigma(M_{k\text{-AP}})=2$ in order to obtain the desired volume.
\end{proof}

\subsection{$B_h[g]$-sets}\label{ssec:Bh}

A polytope with unimodular matrix (namely, each quadrangular submatrix has determinant either $0$ or $\pm 1$) is integral~\cite{Sch86}. We start proving that the polytope associated to $B_h[g]$-sets is integral, hence we can use the usual interpolation technique in polynomials. Note that it is easy to verify that the matrix associated to $B_h[g]$-sets is strictly balanced.

\begin{proposition}\label{prop: unimodularity}
The polytope associated to $B_h[g]$-sets is integral.
\end{proposition}

\begin{proof}
Let $M_{B_h[g]}$ be defined as in Equation~\eqref{eq: matrix_Bhg}. We recall Equation~\eqref{eq:polytopetomatrix} and note that the polytope $P_{M_{B_h[g]}}$ can be written as
\begin{equation*}
\mP_{M_{B_h[g]}}=\{\bx: M_{B_h[g]} \cdot \bx \leq \textbf{0}\}\cap\{\bx: (-M_{B_h[g]})\cdot \bx \leq \textbf{0}\}\cap[0,1]^m\subset \mathbb{R}^{m}.
\end{equation*}
and hence has matrix representation form $\{\textbf{x}\in \mathbb{R}^k: P \cdot \textbf{x} \geq \textbf{b}\}$ where
\begin{equation}\label{eq:system-Bhg}
P = \left( \begin{array}{c} M \Bstrut \\ -M \Tstrut \Bstrut \\ I_{h(g+1)} \Tstrut \Bstrut  \\ -I_{h(g+1)} \Tstrut	\end{array} \right) \quad \text{and} \quad  \textbf{b} = (0, \overset{2g}{\ldots}, 0, 0, \overset{h(g+1)}{\ldots}, 0, -1, \overset{h(g+1)}{\ldots}, -1)^T.
\end{equation}
It follows that we only need to prove that all minors of the matrix belong to the set $\{0,\pm 1\}$. Observe that we can reduce our argument to minors with entries in the topmost part of the matrix (namely the matrix $M$). We argue by induction on the size of the minor. The result is clear for minors of size $1$, as the entries of the matrix belong to $\{0,\pm 1\}$. Assume that the result is true for every minor of size at most $k$, and let us show that the result is also true for $k$. With this purpose we use the fact that every column of $M$ has at most two elements different from $0$.

Consider the first row of the minor under study. If all elements are equal to $0$, the minor is equal to $0$. If there exist a unique element different from $0$, we apply induction by developing the determinant along the row. Finally, let us assume that there exist in the first row at least two elements different from 0. Finally, observe that:

\begin{itemize}
\item [1.] if these two elements in the first row are equal the corresponding columns are linearly dependent, and the determinant is equal to $0$.
\item [2.] if these two elements are different, the column where $1$ belongs just contain $0$ since we are considering the first row: by construction a minor cannot have a $-1$ below a $1$ entry. Hence we can develop the determinant by this column and we apply induction.
\end{itemize}
With this analysis we cover all possible cases, and the proof is finished.
\end{proof}

We continue computing the number of solutions of the system of equations $M_{B_{h}[g]}\cdot \bx=0$ such that the coordinates of $\bx$ belong to $[n]\cup \{0\}$ by means of the inclusion-exclusion method. For this, fix some $j$ coordinates and index them by the set $J \subseteq [h]$. If these coordinates are strictly greater than $n$, we can rewrite our equation $x_1 + \ldots + x_h = k$ as $\sum_{i \in J} (x_i-n)+ \sum_{i \in [h] \backslash J} x_i = y_1 + \ldots + y_h = k - (n+1)j$ where the coordinates of $\by$ still belong to $[n]\cup \{0\}$. It follows that the number of solutions of $x_1 + \ldots + x_h = k$ such that at least the $j$ fixed coordinates are strictly greater than $n$ is equal to
$$\binom{k-(n+1)j+h-1}{h-1}.$$
Now we can apply the inclusion-exclusion method. For this write $k = k_1 n + k_2$ where $0 \leq k_1 \leq h-1$ and $1 \leq k_2 \leq n$ or $k_1=0,\, k_2=0$ (so they are uniquely determined for each $k$). It follows that the number of solutions of the equation $x_1 + \ldots + x_h = k$ with $x_i \in [n]\cup \{0\} $ is equal to
$$a(k_1) = \sum_{j=0}^{k_1} (-1)^{j}\binom{h}{j} \binom{(k_1-j)n+k_2-j+h-1}{h-1}.$$
Observe that $k$ is smaller or equal than $hn$. Consequently, the total number of integer points in the polytope defined by the equations $M_{B_{h}[g]}\cdot \bx=0$ and each coordinate of $\bx$ belonging to $[n]\cup \{0\}$ is equal to a function
$$f_{h,g}(n)=1+\sum_{k_1=0}^{h-1}\sum_{k_2=1}^{n} \left( a(k_1) \right)_{g+1},$$
where we have used the Pochhammer symbol $\left( a(k_1) \right)_{g+1} = a(k_1) \left( a(k_1)-1 \right) \cdots \left( a(k_1)-g \right)$. Now the argument used in the case of $k$-AP does not work, as expressions are more involved. However, we can apply an interpolation argument to obtain the dominant term of $f_{h,g}(n)$:  by Proposition~\ref{prop: unimodularity} and Theorem~\ref{thm:ethm}, $f_{h,g}(n)$ is a polynomial of degree $d = (h-1)(g+1)+1$ with coefficients $a_0$, $a_1$ \ldots, $a_d$. Hence, the values $f_{h,g}(0), f_{h,g}(1),\dots, f_{h,g}(d-1)$ determine $f_{h,g}(n)$ through the Vandermonde-matrix
\begin{equation}\label{eq:det}
\left(
  \begin{array}{cccc}
    1 & 0 & \cdots & 0 \\
    1 & 1 & \cdots & 1 \\
    1 & 2 & \cdots & 2^{d-1} \\
    \vdots &  \vdots & \cdots & \vdots \\
    1 & d-1 & \cdots & (d-1)^{d-1} \\
  \end{array} \right) \cdot \left(
  \begin{array}{c}
    a_0 \\ a_1 \\ a_2 \\ \vdots \\ a_d
  \end{array} \right) = \left(
  \begin{array}{c}
    f_{h,g}(0) \\ f_{h,g}(1) \\ f_{h,g}(2) \\ \vdots \\ f_{h,g}(d-1)
  \end{array} \right)
\end{equation}
Note that by Equation~\eqref{eq:leadingcoefficient} we have $a_d = \Vol{\mP_{M_{B_h[g]}}}$. This coefficient can easily be determined using for example Cramer's rule. Actual values of the volume of the polytope with $h,g\leq 6$ are computed in Table~\ref{table:vol-Bh[g]} using this equation.

\SP

We recall that a detailed study for the threshold in $B_h[1]$ sets can be found in~\cite{GJLR99}. In this work Godbole et al. studied the random variable that counts the number of solutions $(\textbf{a},\textbf{b})=(a_1,\dots,a_h, b_1, \dots,b_h)$ of the equation
\begin{equation}\label{janson_sol}a_1+a_2+\dots+a_h=b_1+b_2+\dots+b_h,\end{equation}
with $a_1\leq a_2\leq \dots\leq a_h$, $b_1\leq b_2\leq \dots\leq b_h$ and $\textbf{a}<\textbf{b}$ with respect to the lexicographic order. They obtained the volume of the associated polytope by means of trigonometric sums and Fourier analytic methods. More precisely, this volume is given by Equation~$(16)$ in~\cite{GJLR99}:
$$\kappa_h=\frac{1}{2(h!)^2 (2h-1)!}\sum_{j=0}^{h-1}(-1)^j\binom{2h}{j}(h-j)^{2h-1}.$$
As we have seen before, it suffices to study carefully the number of proper solutions. Therefore, in terms of our approach, this result can be translated into
$$\Vol{\mP_{B_h[1]}} = \sigma(M_{B_h[1]})\kappa_h = 2(h!)^2\kappa_h = \frac{\sum_{j=0}^{h-1}(-1)^j\binom{2h}{j}(h-j)^{2h-1}}{(2h-1)!},$$
since for every ordered solution to~\eqref{janson_sol} we must count $2(h!)^2$ different solutions (obtained by permuting the $a_i$ and the $b_j$ coordinates, and then considering the symmetric solution $(\textbf{b},\textbf{a})$). These constants correspond to the first column in the following table ($g=1$). Closed formulas for bigger values of $g$ seem to be much more involved.

\begin{table}[htb]
\begin{center}
  \begin{tabular}{c|ccccc}
  $h\diagdown g$ &1&2&3&4&5 \Bstrut \\
  \hline

  2 &$\frac{2}{3}$&  $\frac{1}{2}$& $\frac{2}{5}$ & $\frac{1}{3}$ & $\frac{2}{7}$ \Tstrut \Bstrut \\
  3 & $\frac{11}{20}$ &  $\frac{12}{35}$&  $\frac{379}{1680}$ &$\frac{565}{3696}$ & $\frac{6759}{64064}$ \Tstrut \Bstrut \\
  4 & $\frac{151}{315}$  & $\frac{1979}{7560}$ & $\frac{40853}{270270}$ & $\frac{200267}{2223936}$ & $\frac{825643615}{15084957888}$ \Tstrut \Bstrut \\
  5 & $\frac{15619}{36288}$ & $\frac{4393189}{20756736}$ & $\frac{1865002207}{16937496576}$ & $\frac{342366164065}{5792623828992}$ & $\frac{689860777579903}{21316855690690560}$ \Tstrut \Bstrut \\
  6 & $\frac{655177}{1663200}$ & $\frac{45515121}{256256000}$  & $\frac{1549892743123}{18284797440000}$ & $\frac{1931111804640401}{46260537523200000}$ & $\frac{31400953991819767493}{1497176036400844800000}$ \Tstrut \Bstrut \\
  \end{tabular}\bigskip
\caption{Volumes for different families of $B_{h}[g]$ sets.}\label{table:vol-Bh[g]}
\end{center}
\end{table}

\section{Related questions}\label{related}

The problem considered in this paper could be rephrased in a more general setting. Let $\mQ$ be an infinite sequence of integers. Let $\mA$ be a random set in $[n]$, and $M \in \mM_{r,m}(\ZZ)$ an admissible matrix. Does there exist a threshold function for the property ``\emph{$\mA^m$ contains a non-trivial solution $\bx$ with $M\cdot \bx \in \mQ^{r}$ }''? Observe that this paper has dealt with the case $\mQ=\{0\}$. It is clear that we need extra assumptions on the the matrix $M$: for instance, the system of equations with matrix
$$M=\left(
    \begin{array}{cc}
      2 & -4 \\
      1 & 2 \\
    \end{array}
  \right)
$$
is admissible, but $M\cdot  \bx \in \mQ^2$ when $\mQ=2\mathbb{N}+1$ is not possible. The problem of characterizing those matrices which are admissible for a given sequence $\mQ$ or, on the contrary, characterizing those sequences that are \emph{admissible for a fixed system} is a problem far from being trivial. Nevertheless, even for very simple systems (such as $x_1-x_2$) and well studied sequences (like the squares or the primes) the study of large sets which avoid this condition is very involved.

For example, S\'{a}rk\"{o}zy~\cite{Sa78} showed that every set with positive upper density contains at least two elements whose difference is a square, see also~\cite{Ly13}. It is, in fact, conjectured that for every $\epsilon>0$ there exists a set in $[n]$ whose differences are never a square and has size $n^{1-\epsilon}$. Ruzsa~\cite{Rz84} proved this conjecture for every $\epsilon\geq 0.267$.

In the presented approach, however, some things can be said. For example, consider the equation $x_1-x_2$ and the sequence of $k$-th powers $\mQ=\{x^k: x\in \mathbb N\}$ (the same arguments could be applied to more general sequences, like prime numbers or powers of $2$ among others).

Then, it is obvious that, if we denote by $S_{\mQ}(n)=\{\bx=(x_1,x_2)\in [n]^2: x_1-x_2\in\mQ\}$ the set of solutions,
$$|S_{\mQ}(n)|=\sum_{q\in\mQ(n)}(n-q)= n|\mQ(n)|-\sum_{q\in \mQ(n)}q=\int_{0}^n x^{1/k}dx= \frac{k}{k+1}n^{1+1/k}(1+o(1))$$
by Abel's summation formula.

It is easy to see that if $\mA$ a random set of $[n]$, where every element is chosen uniformly at random with probability $p$, then $p=n^{-(k+1)/(2k)}$ is a threshold function for the property ``$x_1-x_2\in \mQ$''. The proof follows from the same ideas of Theorem~\ref{thm:threshold}. Once again, if we denote by $E_{\bx}$ the event $\bx\in\mA^2$ and $\II_{\bx}$ is the associated indicator random variable, it is clear that the expected value for the random variable
$$\bX=\sum_{\bx\in S_{\mQ}(n)}\II_{\bx},$$
is  $O\left(n^{(k+1)/{k}}p^2\right)$. Hence, taking $p=o\left(n^{-(k+1)/(2k)}\right)$ this expected value tends to $0$.

For the second part, we observe that
$$\Delta=O\left(n|\mQ(n)|^2 p^3\right)=O(n^{\frac{k+2}{k}}p^3)$$
and therefore taking $p\gg n^{-\frac{k+1}{2k}}$ we obtain that $\Delta=o\left(\mathbb{E}[\bX]^2\right)$. Consequently, $\bX \sim \mathbb{E}[\bX]$ asymptotically almost surely.

The methodology developed to deal with systems of linear equations could be adapted to treat similar problems in other directions. The same arguments could be adapted in the context of  finite fields: despite the extra conditions we need to demand to the system (in order to have maximum rank), we do not need an Ehrhart's type result in this context.

\SP

\paragraph{\textbf{Acknowledgments}: The authors thank Christian Elsholtz for providing many references concerning linear systems of equations. We also thank Katy Beeler, Arnau Padrol and Lluis Vena for valuable input and fruitful discussions as well as Javier Cilleruelo for a detailed reading of the manuscript and support. Furthermore, we are thankful to the anonymous referees for detailed feedback and constructive suggestions. We would like to thank the \emph{Combinatorics and Graph Theory} group at the Freie Universität Berlin where part of this research took place for working conditions and hospitality.}

\bibliography{bib}
\bibliographystyle{abbrv}

\end{document}

\cite{LRS10}